\documentclass[reqno,A4paper]{amsart}

\usepackage{amsmath}
\usepackage{amssymb}
\usepackage{amsthm}
\usepackage{enumerate}
\usepackage{mathrsfs} 
\usepackage{eqlist}
\usepackage{array}
\usepackage{xcolor}
\newtheorem{Theo}{Theorem}[section]
\newtheorem{Lem}{Lemma}[section]

\newtheorem{Assumpt}{Assumption}[section]

\theoremstyle{definition}

\theoremstyle{Rmq}
\newtheorem{Rmq}{Remark}[section]

\theoremstyle{plain}

\theoremstyle{definition}

\numberwithin{equation}{section}

\begin{document}

\title[Nonparametric estimation of SIR by  the $ k$-NN kernel method]%
{Nonparametric estimation of sliced inverse regression by the   $ k$-nearest neighbors kernel method}
\author[L.  Bengono Mintogo \and E.D.D. Nkou \and G.M. Nkiet]%
{Luran  Bengono Mintogo\and  Emmanuel de Dieu  Nkou  \and Guy Martial Nkiet}

\newcommand{\acr}{\newline\indent}

\address{\llap{ }Laboratoire de Probabilit\'es, Statistique et Informatique \acr
                  Unit\'e de Recherche en Math\'ematiques et Informatique\acr
                  Université des Sciences et Techniques de Masuku\acr
                  BP 943  Franceville\acr
                  GABON}
\email{luranbengono@gmail.com, emmanueldedieunkou@gmail.com, guymartial.nkiet@univ-masuku.org }




\begin{abstract}
We investigate nonparametric estimation of  sliced inverse regression  (SIR) via the $k$-nearest neighbors approach with a kernel. An estimator of the covariance matrix of  the conditional expectation of the explanatory random vector given the response is then introduced, thereby allowing to estimate the effective dimension reduction (EDR) space.  Consistency of the proposed estimators is proved through derivation of asymptotic normality. A simulation study, made in order to  assess the finite-sample behaviour of the proposed method and to compare it to the kernel estimate, is presented. 

\bigskip

\noindent R{\tiny \'ESUM\'E}. Nous abordons l'estimation non-param\'etrique de la r\'egression inverse par tranche (SIR) par la m\'ethode des $k$ plus proches voisins \`a noyau. Un estimateur de la matrice des covariances de l'esp\'erance conditionnelle  du vecteur al\'eatoire explicatif   sachant  la variable al\'eatoire \`a expliquer  est alors  introduit, permettant ainsi d'estimer l'espace de r\'eduction effective de la dimension (EDR). La convergence des estimateurs propos\'es est prouv\'ee \`a travers la d\'etermination de leur normalit\'e asymptotique.  Une \'etude par simulations,   permettant d'évaluer   les performances  \`a taille d'\'echantillon finie  de la m\'ethode propos\'ee et de la comparer \`a l'estimation par noyau, est pr\'esent\'ee.
\end{abstract}

\maketitle
\section{Introduction}
\noindent Dimension reduction has become a major concern  in statistical multivariate analysis since the introduction of sliced inverse regression (SIR)  by \cite{li91}. It is an approach allowing to handle a regression with a reasonable number of explanatory variables  by projecting an  initial multidimensional  predictor onto a subspace of lower dimension than that of this predictor. For doing that, \cite{li91} introduced the model
\begin{equation}\label{mod}
Y=F(\beta_1^{\top}X,...,\beta_{N}^{\top}X,\varepsilon),
\end{equation}
where  $Y$ is a scalar response,  $X$ is a $d$-dimensional random vector containing the initial explanatory variables as coordinates, $N$ is an integer satisfying  $1\leqslant N<d$, $\beta_1,\ldots,\beta_N$ are unknown vectors in $\mathbb{R}^d$,  $\varepsilon$ is a   scalar random variable  independent  of $X$, $u^\top$ denotes the transpose of the vector $u$,  and $F$ is an arbitrary unknown function. This model allows to perform the aforementioned  reduction of the number of predictors since $N<d$, but   requires an appropriate  choice of the vectors $\beta_1,\ldots,\beta_N$ so that the projection onto the subspace  they span, called the effective dimension reduction (EDR) space, retains a maximum of information on $Y$.   Estimating   the EDR space  is then the most crucial issue related to this model and has  been tackled in several works. Li\cite{li91} proposed  an approach  based on estimating the covariance matrix of the conditional expectation $\mathbb{E}\left(X\vert Y\right)$ from slicing the range of $Y$.  Although there exist several  alternative methods,  it  stills the most recognized method in sufficient dimension reduction, and  has seen numerous   developments such as asymptotic studies (\cite{hsing92,zhu95}), extensions to the multiple scalar responses case (\cite{coudret14,li03}), estimation of the dimensionality $N$ (\cite{ferre98,nkiet08,schott94,velilla98}), extension to the high-dimensional framework (\cite{zhu06}). Among the alternative methods, there is the sliced average variance estimation (SAVE) method of \cite{cook91}, the parametric inverse regression (PIR) method of \cite{bura01}, the central mean subspace (CMS) estimation method introduced in \cite{cook02}.

 Based on the fact that the aforementioned covariance  matrix can be expressed  as a function of the density of $Y$ and regression functions,    \cite{zhu96} proposed a nonparametric approach for estimating the EDR space  by using Nadaraya-Watson type estimators. This enlighted for the first time the possibility and the interest of using nonparametric estimation methods for estimating SIR. Later, \cite{nkou19} established the strong consistency of this approach, \cite{nkou22} proposed another nonparametric estimation approach based on wavelets and \cite{nkou23} used recursive kernel estimators. However, the practical choice of the bandwidth on which the estimators used in \cite{zhu96}  rely
is not straightforward and stills a challenging issue. So, there is an interest in using  alternative estimators  which do not require to make  such a choice. Among them, the $k$-nearest neighbors ($k$-NN) kernel estimators have attracted particular attention. They have the same form than the kernel estimators, but with bandwidth replaced by the Euclidean distance between the point to which the estimator is calculated and the $k$th nearest neighbor  of this point among the observations. Earlier works on these estimators  go back to \cite{moore77} for density estimation and to \cite{collomb79} for the case of regression function.  It is true that the practical implementation of these estimators requires choosing the number of neighbors, but this choice is easier to make since it is to be made on a finite  set of integers. To the best of our knowledge, these estimators have never been used for performing dimension reduction methods.

In this paper,    we propose  an estimation method for the EDR space related to the  model \eqref{mod}, based on $k$-NN kernel   estimates of the density and regression functions involved in the aforementioned covariance matrix. 

The rest of the paper is organized as follows.  In Section \ref{knn},  we define the tackled estimators   based on $k$-NN kernel method. Their asymptotic properties    are then given in   Section \ref{asymp}.  Section \ref{simul} is devoted to a simulation study made in order to evaluate the performance of the proposal with comparison to the  method of \cite{zhu96}.   The proofs of the main results are postponed in Section \ref{proofs}. 
 
\section{Estimation based on the $ k$-NN method}\label{knn}

\subsection{The matrix to be estimated}
The most  crucial issue in handling the model \eqref{mod} is to estimate the EDR space which is known  to be spanned, under some conditions,  by eigenvectors of a matrix that will now be specified. We assume that the explanatory random vector 
\[
X=\left(X_1,\ldots,X_d\right)^\top,
\] 
satisfies  $\mathbb{E}\left(\Vert X\Vert^2\right)<\infty$, where $\Vert\cdot\Vert$ denotes the usual Euclidean norm of $\mathbb{R}^d$, and, without loss of generality,   $\mathbb{E}\left(X\right)=0$. Then  it is known  from \cite{li91}  that, under some conditions,  the EDR space is spanned by  eigenvectors associated to the  largest eigenvalues of the covariance matrix $\Lambda$ of the conditional expectation of $X$ given $Y$ denoted, as usual, by $\mathbb{E}\left(X\vert Y\right)$; that is 
\[
\Lambda=\mathbb{E}\left(R(Y)R(Y)^\top\right), 
\]
where   $R(Y)=\mathbb{E}\left(X\vert Y\right)$. Since 
\[
R(Y)=\left(R_1(Y),\ldots,R_d(Y)\right)^\top,
\] 
where $R_\ell(Y)=\mathbb{E}\left(X_\ell\vert Y\right)$ for $\ell\in\{1,\ldots,d\}$, it is clear that 
\[
\Lambda=\Big( \lambda_{\ell j}\Big)_{1\leqslant\ell , j\leqslant d}\,\,\,\textrm{ with }\lambda_{\ell j}=R_\ell(Y)R_j(Y).
\]
Estimation of the EDR space  is, therefore,  reduced to that of $\Lambda$, what is obtained from estimating the $\lambda_{\ell j}$'s. For doing that, \cite{li91} introduced a method, called sliced inverse regression (SIR), based on slicing the range of $Y$ and then estimating $\Lambda$ by a scatter matrix computed by using the dispersion of the observations within the different slices. As the $R_\ell$'s are  in fact   regression functions,  they can be estimated from nonparametric approaches. More specifically, assuming  that each pair $(X_\ell,Y)$ admits a density that we denote by $f_{ (X_\ell,Y)}$, and that the density    $f$    of $Y$ is such that     $f(y)>0$ for all $y\in\mathbb{R}$, we have
\[
R_\ell(y)=\mathbb{E}(X_\ell\vert Y=y)=\frac{g_\ell(y)}{f(y)}
\,\,\textrm{ where }\,\,g_\ell(y)=\int_\mathbb{R} x f_{_{(X_\ell,Y)}}(x,y)\,\,dx.
\] 
Based on this, \cite{zhu96} proposed a nonparametric estimation of $\Lambda$ based on  a Nadaraya-Watson type estimator of $R_\ell$ whereas \cite{nkou22} uses estimation based on  wavelets. In this paper, we use rather estimation based on the $k$ nearest neighbors method with a kernel.
\subsection{The $k$-NN kernel estimator}
 Considering an i.i.d. sample  $\{(X^{(i)},Y_i)\}_{1\leqslant i\leqslant n}$   of     $(X,Y )$,   and putting
$
X^{(i)}=\left(X_{i1},\ldots,X_{id}\right)^\top,
$
we make use of the $k$-NN kernel estimators    of the density  $f$ and the   functions  $g_\ell$, as defined in \cite{bengono24,collomb79,moore77}. More specifically, given a kernel $K\,:\,\mathbb{R}\rightarrow\mathbb{R}$, we consider the estimators $\widehat{f}_{n}$ and $\widehat{g}_{\ell,n}$ of $f$ and $g_\ell$, respectively, given by
\begin{equation*}\label{estf}
\widehat{f}_{n}(y)=\frac{1}{nH_{n}(y)}\sum\limits_{i=1}^{n}K\left(\frac{Y_{i}-y}{H_{n}(y)}\right), 
\end{equation*}
and
\begin{equation*}\label{estgl}
\widehat{g}_{\ell,n}(y)=\frac{1}{nH_{n}(y)}\sum\limits_{i=1}^{n}X_{i\ell}K\left(\frac{Y_{i}-y}{H_{n}(y)}\right),
\end{equation*}
where
\begin{equation*}
H_{n}(y)=\min\left\{h\in\mathbb{R}_{+}^{*}\bigg{/}\sum\limits_{i=1}^{n}\mathbf{1}_{]y-h,y+h[}(Y_{i})=k_{n}\right\},
\end{equation*}
$\ell\in\{1,\ldots,d\}$, and  $\left(k_n\right)_{n\in\mathbb{N}^\ast}$ is a sequence of integers such that $k_n\rightarrow \infty$ as $n\rightarrow \infty$.
 As it was already done in \cite{bengono24,nkou23,nkou19,nkou22,zhu96},  in order to avoid  small values in the denominator, we  consider
\[
f_{b_n}(y)=\max\big\{f(y),b_n\big\}\,\,\textrm{ and }\,\, \widehat{f}_{b_n}(y)=\max\big\{\widehat{f}_n(y),b_n\big\},
\]
where  $\left(b_n\right)_{n\in\mathbb{N}^\ast}$ is a   sequence of positive  real numbers  such that   $\lim_{n\rightarrow \infty}(b_n)=0$, and we estimate the ratio    
\begin{equation}\label{rbnl}
R_{b_n,\ell}(y)=\frac{ g_\ell(y)}{f_{b_n}(y)}
\end{equation}
 by  
\begin{equation*}\label{hatrbj}
 \widehat{R}_{b_n,\ell}(y)= \frac{ \widehat{g}_{\ell,n}(y)}{\widehat{f}_{b_n}(y)}.
\end{equation*}
Then,  putting
\begin{equation*}\label{hatrb}
\widehat{R}_{b_n}(y)=\left(\widehat{R}_{b_n,1}(y),\ldots,\widehat{R}_{b_n,d}(y)\right)^\top,
\end{equation*}
we take as estimator of   $\Lambda$ the random matrix
\begin{equation*}\label{hatlambda}
\widehat{\Lambda}_n=\frac{1}{n}\sum_{i=1}^{n}\widehat{R}_{b_n}(Y_i)\widehat{R}_{b_n}(Y_i)^\top.
\end{equation*}
An estimate of the EDR space is obtained from the spectral analysis of this matrix. Indeed, if  $\widehat{\beta}_\ell$ is   as an eigenvector of $\widehat{\Lambda}_n$ associated with the $\ell$-th largest eigenvalue, we estimate the EDR space by the subspace of $\mathbb{R}^d$  spanned by   $\widehat{\beta}_1,\cdots, \widehat{\beta}_N$.

\bigskip

\begin{Rmq}
The main difference between our approach and that of \cite{zhu96} is about the bandwith of the used estimators. A fixed real bandwith, that only depends on $n$, is used in \cite{zhu96} whereas our approach consists in considering a random bandwith that depends on the observations and the point of estimation. Such estimators, that belong to the broader class of variable kernel estimators, have some advantages over fixed bandwith estimators. In particular, they do not require to make  a prior choice of the bandwith as it is the case for the fixed bandwith estimators. It is true that the choice of the number of neighbors must be made, but this is done in a finite set of integers, which makes the research easier.
\end{Rmq}

\section{Asymptotic properties}\label{asymp}
\noindent  In this section, we first introduce the assumptions needed to obtain the main results of the paper,  and then we state
 theorems that give asymptotic normality for the considered estimators.
\subsection{Assumptions}
We make the following assumptions:

\bigskip
\begin{Assumpt}\label{as0}
	$\mathbb{E}\left(\left\Vert X\right\Vert^4\right)<\infty$.  
\end{Assumpt}
\begin{Assumpt}\label{as7}
	 There exists a sequence $M_{n}$ of strictly positive numbers such that $M_{n}\sim\sqrt{\log(n)}$ and, for $\ell\in\{1,\ldots,d\}$,  $\max_{1\leqslant  i\leqslant n}\left\vert X_{i\ell}\right\vert\leqslant M_{n}$.
\end{Assumpt}
\begin{Assumpt}\label{as1}
	The density $f$ of $Y$ is bounded from below: there exists $c_{0}>0$ such that $\inf_{y\in\mathbb{R}}f(y)\geqslant c_{0}.$
\end{Assumpt}

\begin{Assumpt}\label{as2}
	The density $f$ is bounded and belongs to the class $\mathscr{C}(c,3)$ of functions $\varphi:\mathbb{R}\rightarrow\mathbb{R}$ that are 3 times differentiable with third derivative satisfying the Lipschitz condition:
$$\left|\varphi^{(3)}(y+u)-\varphi^{(3)}(y)\right|\leqslant c\left|u\right|,$$
where $c>0$.
\end{Assumpt}
\begin{Assumpt}\label{as3}
	 For $\ell\in\left\{1,\ldots,d\right\}$, the functions $g_{1,\ell}$ and $g_{2,\ell}$ defined by  
\begin{equation}\label{g1lg2l}
g_{1,\ell}(y)=\mathbb{E}\left(X_{\ell}\mathbf{1}_{\{X_{\ell}\geqslant 0\}}|Y=y\right)f(y),\,\,\,g_{2,\ell}(y)=\mathbb{E}\left(-X_{\ell}\mathbf{1}_{\{X_{\ell}<0\}}|Y=y\right)f(y),
\end{equation} 
are bounded and belong to the class $\mathscr{C}(c,3)$ previously defined.
\end{Assumpt}
\begin{Assumpt}\label{as8}
	For all $\ell \in \{1, \ldots, d\}$ and all $m\in\{1,2\}$, the function $R_{m,\ell}$ defined by
\begin{equation}\label{r1lr2l}
R_{m,\ell}(y)=\frac{g_{m,\ell}(y)}{f(y)}
\end{equation}
satisfies:
\begin{enumerate} 
\item[(i)]  $\mathbb{E}\left[R_{m,\ell}^4(Y)\right]< \infty$.
\item[(ii)] $\left|R_{m,\ell}(y+u)-R_{m,\ell}(y)\right|\leqslant c\left|u\right|$.
\end{enumerate}
\end{Assumpt}
\begin{Assumpt}\label{as9}  
$\sqrt{n}\,\mathbb{E}\left[\left|R_{m,\ell}(Y)R_{s,j}(Y)\right|\mathbf{1}_{\{f(Y)\leqslant a_n \}}\right] = \mathrm{o}(1)$  for any $(\ell,j)\in\{1,\ldots,  d\}^2$, $(m,s)\in\{1,2\}^2$ and any sequence $(a_n)_{n\in\mathbb{N}^\ast}$ such that $a_n\sim b_n $ as $n\rightarrow \infty$.
\end{Assumpt}
\begin{Assumpt}\label{as4}
The kernel $K:\mathbb{R}\rightarrow\mathbb{R}$ satisfies the following properties:
\begin{enumerate} 
\item[(i)] $K$ is bounded, that is $G=\sup_{t\in\mathbb{R}}\left|K(t)\right|<\infty$.
\item[(ii)] $K$ is symetric with respect to 0, that is $K(t)=K(-t)$, $\forall t\in\mathbb{R}$.
\item[(iii)] $\int_{\mathbb{R}}K(t)\,\,dt=1$.
\item[(iv)] $K$ is of order 3, that is $$\int_{\mathbb{R}}t^{m}K(t)\,\,dt=0,\hspace{0.2cm}m\in\left\{1,2,3\right\}.$$ 
\item[(v)] $$\int_{\mathbb{R}} \left|K(t)\right|\,\,\,\,dt<\infty,\int_{\mathbb{R}} \vert t\vert\,\,\vert K(t)\vert \, dt<\infty\,\,\textrm{ and }\,\,\int_{\mathbb{R}}t^{4}\left|K(t)\right|\,\,dt<\infty.$$
\item[(vi)] $\forall t\in\mathbb{R}$, $\forall a\in \left[0,1\right]$, $K(at)\geqslant K(t)$.
\end{enumerate}
\end{Assumpt}
\begin{Assumpt}\label{as5}
The number $k_{n}$ of neighbors is   such that $k_{n}\sim  n^{c_{1}}$, where $1/2<c_{1}<9/10$.
\end{Assumpt}
\begin{Assumpt}\label{as6}
The sequence $\left(b_{n}\right)_{n\in\mathbb{N}^{*}}$ satisfies $b_{n}\sim n^{-c_{2}}$ with $0<c_{2}<1/10$, $c_1>c_2+3/4$ and $c_1+c_2/4<7/8$.
\end{Assumpt}
\begin{Assumpt}\label{as10}
The eigenvalues $\nu_1,\ldots,\nu_d$ of $\Lambda$ satisfy $\nu_1>\nu_2>\cdots>\nu_d>0$.
\end{Assumpt}

\begin{Rmq}
	Assumption \ref{as0} is a  classical one  in the literature on asymptotic studies in multivariate analysis. In the context of dimension reduction techniques, it was made, for instance,  in \cite{zhu96}.  Assumption \ref{as7} is weaker than Assumption 3.1 of \cite{nkou19} and Assumption 3.1  of \cite{nkou22}  where it is supposed that $X$ is bounded. For instance, it has  been considered in \cite{bengono24,nkou23}.  Assumption \ref{as1} has been introduced in some works on nonparametric estimation such as \cite{bengono24, ebende24,zhu07}. Assumption \ref{as2} just is the first part of Condition 1 of \cite{zhu96} whereas Assumption \ref{as3} implies its  second part. Indeed, since $g_\ell=g_{1,\ell}-g_{2,\ell}$, it follows that
\[
\left\vert g_\ell(y)-g_\ell(z)\right\vert\leqslant \left\vert g_{1,\ell}(y)-g_{1,\ell}(z)\right\vert+\left\vert g_{2,\ell}(y)-g_{2,\ell}(z)\right\vert .
\]
So, if $g_{1,\ell}$ and $g_{2,\ell}$ belong to   $\mathscr{C}(c,3)$, then $g_\ell$ belongs to $\mathscr{C}(2c,3)$. It is the same for Assumptions \ref{as8} and \ref{as9} which imply Assumptions 3 and 7 of \cite{nkou22} and Condition 6 of \cite{zhu96} since we have $R_\ell=R_{1,\ell}-R_{2,\ell}$ which yields
\[
\left\vert R_\ell(y)-R_\ell(z)\right\vert\leqslant \left\vert R_{1,\ell}(y)-R_{1,\ell}(z)\right\vert+\left\vert R_{2,\ell}(y)-R_{2,\ell}(z)\right\vert ,
\]
\[
\left\vert R_\ell(y)\right\vert\leqslant \left\vert R_{1,\ell}(y) \right\vert+\left\vert R_{2,\ell}(y) \right\vert ,
\]
and
\[
 R^4_\ell(y) \leqslant 8\left(  R_{1,\ell}^4(y)  +  R_{2,\ell}^4(y) \right).
\]
So, if $R_{1,\ell}$ and $R_{2,\ell}$ belong to   $\mathscr{C}(c,3)$, then $R_\ell$ belongs to $\mathscr{C}(2c,3)$. Moreover, if $\mathbb{E}\left[R_{m,\ell}^4(Y)\right]< \infty$ for $m\in\{1,2\}$, then $\mathbb{E}\left[R_{\ell}^4(Y)\right]< \infty$, and  if  $\sqrt{n}\,\mathbb{E}\left[\left|R_{m,\ell}(Y)R_{s,j}(Y)\right|\mathbf{1}_{\{f(Y)\leqslant a_n \}}\right] = \mathrm{o}(1)$ for $(m,s)\in\{1,2\}^2$ and  $a_n\sim b_n $, then $\left\vert\sqrt{n}\,\mathbb{E}\left[R_{\ell}(Y)R_{j}(Y)\mathbf{1}_{\{f(Y)\leqslant b_n \}}\right]\right\vert = \mathrm{o}(1)$ because
\begin{align*}
&\left\vert\sqrt{n}\,\mathbb{E}\left[R_{\ell}(Y)R_{j}(Y)\mathbf{1}_{\{f(Y)\leqslant b_n \}}\right]\right\vert \\
&\leqslant \sqrt{n}\,\mathbb{E}\left[\left\vert R_{\ell}(Y)R_{j}(Y)\right\vert\mathbf{1}_{\{f(Y)\leqslant b_n \}}\right]\\
&\leqslant \sqrt{n}\,\mathbb{E}\left[\left\vert R_{1,\ell}(Y)R_{1,j}(Y)\right\vert\mathbf{1}_{\{f(Y)\leqslant b_n \}}\right]+\sqrt{n}\,\mathbb{E}\left[\left\vert R_{1,\ell}(Y)R_{2,j}(Y)\right\vert\mathbf{1}_{\{f(Y)\leqslant b_n \}}\right]\\
&+\sqrt{n}\,\mathbb{E}\left[\left\vert R_{2,\ell}(Y)R_{1,j}(Y)\right\vert\mathbf{1}_{\{f(Y)\leqslant b_n \}}\right]+\sqrt{n}\,\mathbb{E}\left[\left\vert R_{2,\ell}(Y)R_{2,j}(Y)\right\vert\mathbf{1}_{\{f(Y)\leqslant b_n \}}\right].
\end{align*}
The conditions $(i)$ to $(v)$ of Assumption   \ref{as4} are classical ones and was considered, for instance, in \cite{bengono24,ebende24,nkou19,zhu96,zhu07}. The condition $(vi)$ often arises in the literature on $k$-NN kernel estimators. Assumptions \ref{as5} and \ref{as6} are of types that have already be considered in the litterature on nonparametric estimation in dimension reduction methods (e.g., \cite{bengono24,ebende24,nkou19,nkou22,zhu96,zhu07}). 
\end{Rmq}
\subsection{Results}
Here, we give the main results of the paper which give asymptotic normality for the estimator $\widehat{\Lambda}_n$ introduced previously  and also, as a consequence, that of its eigenvectors under some specified conditions.

\bigskip

\begin{Theo}\label{theo1}
Under the  assumptions \ref{as0} to  \ref{as6}, we have
\[
\sqrt{n}\left(\widehat{\Lambda}_n-\Lambda\right)\stackrel{\mathscr{D}}{\rightarrow}\mathcal{H}
\]
as $n\rightarrow \infty$, where $\stackrel{\mathscr{D}}{\rightarrow}$ denotes convergence in distribution, $\mathcal{H}$ is a random variable having a normal distribution  in the space $\mathscr{M}_d(\mathbb{R})$ of $d\times d$ matrices, such that, for any $A=\left(a_{\ell j}\right)\in \mathscr{M}_d(\mathbb{R})-\{0\}$,    one has   $\mathrm{Tr}\left(A^\top\mathcal{H}\right) \leadsto \mathcal{N}(0,\sigma^2_{A})$ with:  
\begin{equation}\label{siga}
\sigma_A^2 =Var\left( \sum_{\ell=1}^d\sum_{j=1}^d\frac{a_{\ell j}}{2}\left(X_{\ell }R_j(Y)+X_{j}R_{\ell}(Y)\right)\right).
\end{equation}
\end{Theo}

\bigskip

\noindent From Theorem \ref{theo2}, we can derive the asymptotic normality of the   eigenvectors. For $(\ell,j)\in\{1,\cdots,d\}^2$, we put  $\beta_\ell=\left(\beta_{\ell 1},\ldots,\beta_{\ell d}\right)^\top$ and  we consider the random variable
\begin{equation*}\label{w}
\mathcal{W}_{\ell j}=\left(\sum_{\stackrel{r=1}{r\neq j}}^d\frac{\beta_{rj}}{\nu_\ell-\nu_r}\right)\,\sum_{p=1}^d\sum_{q=1}^d\frac{\beta_{\ell p}\beta_{\ell q}}{2}\Big(X_{q}R_p\left(Y\right)+X_p R_{q}\left(Y\right)\Big).
\end{equation*}
Then, we have:
\bigskip

\begin{Theo}\label{theo2}
Under the  assumptions \ref{as0} to  \ref{as10},  we have for any $\ell\in\{1,\ldots,d\}$,
$
\sqrt{n}\left(\widehat{\beta}_\ell-\beta_\ell\right)\stackrel{\mathscr{D}}{\rightarrow}\mathcal{N}(0,\Sigma_\ell)
$  as $n\rightarrow \infty$,
where $\Sigma_\ell$ is the $d\times d$ covariance matrix of the random vector $ \mathcal{W}_{\ell}=\left(\mathcal{W}_{\ell 1},\cdots,\mathcal{W}_{\ell d}\right)^\top$.
\end{Theo}

\section{Simulation results}\label{simul}
\noindent In order to observe the performance of the introduced  method for estimating SIR, and to compare it with  the method of \cite{zhu96} based on kernel estimates, we made simulations within   frameworks  corresponding to the following models with dimension $d=5$:

\bigskip

\noindent\textbf{Model 1:} $Y=X_{1}+X_{2}+X_{3}+X_{4}+\varepsilon$;

\medskip

\noindent\textbf{Model 2:} $Y=X_{1}\left(X_{1}+X_{2}+1\right)+\varepsilon$;

\medskip

\noindent\textbf{Model 3:}  $Y=X_{1}\left(0.5+\left(X_{2}+1.5\right)^{2}\right)^{-1}+\varepsilon$.

\bigskip

\noindent In these models, which come from \cite{nkou22}, the predictor $X=(X_1,\ldots,X_5)^\top$ is generated from a multivariate normal distribution $\mathcal{N}\left(0,I_{5}\right)$, where $I_{5}$ is the $5\times 5$  identity matrix, $\varepsilon$ is generated independently of $X$ from a standard normal distribution, and $Y$ is computed according to these models. In Model 1, we have $N=1$, $\beta_1=\left(1,1,1,1,0\right)^\top$ whereas Model 2 and Model 3 correspond to the case of $N=2$, $\beta_1=\left(1,0,0,0,0\right)^\top$, and respectively: 
\[
\beta_2=\left(1,1,0,0,0\right)^\top,\,\,\,\beta_2=\left(0,1,0,0,0\right)^\top.
\]
\noindent We generated $500$ independent samples of size $n=50,100,200,400$ from the above models.  For each of these samples we computed  estimates $\widehat{\beta}_1$ and $\widehat{\beta}_2$  of $\beta_1$ and $\beta_2$, and then the distance $\mathcal{D}$  between the true  EDR space and its estimation given by  
\[
\mathcal{D}^2=\textrm{Tr}\bigg[\left(P-\widehat{P}\right)^2\bigg],
\] 
where   $P$  (resp. $\widehat{P}$) denotes the projector onto the space spanned by $\beta_1$  (resp.  $\widehat{\beta}_1$)  in case of  Model 1, and by $\{\beta_1,\beta_2\}$ (resp.  $\{\widehat{\beta}_1,\widehat{\beta}_2\}$) in case of Model 2 or Model 3. Then, the average and standard deviation of $\mathcal{D}$  over the   $500$   replicates are computed   in order to  assess the performance of the  methods. The two methods were used with $b_n= n^{-0.09}$ and  different kernels : Gaussian kernel, Epanechnikov kernel, biweight kernel, triweight kernel and triangular kernel. For our method, we took $k_n$ as the integer part of $n^{0.85}$. The obtained results are reported in Table \ref{tab1} and Table \ref{tab2}. It can be seen that  both methods give equivalent results. Indeed,   when one outperforms  the other, the difference is generally very small. This shows that our method can be seen as a good alternative to the classical method based on kernel estimates.

\begin{table}[htbp]
\centering 
\caption{Average and standard deviation of $\mathcal{D}$ over $500$  replicates, with sample size $n=50,100$.}
	\label{tab1} 
\medskip
	\begin{tabular}{cccccccccc}
\hline\hline
& & & & & & & & &\\
Sample size & Kernel& \multicolumn{2}{c}{Model 1}& & \multicolumn{2}{c}{Model 2}& & \multicolumn{2}{c}{Model 3}\\
\cline{3-4}\cline{6-7}\cline{9-10}
& & $k$-NN& Kernel& &$k$-NN &Kernel & & $k$-NN&Kernel\\
\hline
& & & & & & & & &\\
 &Gaussian&  0.557&0.462& &1.232   & 1.166& &1.240& 1.248\\
& &(0.202) &(0.166) & &(0.208) &(0.212) & &(0.194) &(0.178)\\
& & & & & & & & &\\
 & Epanechnikov 	& 0.504 &0.534 &  &1.232  &1.188  & & 1.238  &  1.236\\
& &(0.180) &(0.200) & &(0.188) &(0.212) & &(0.202) &(0.194)\\
& & & & & & & & &\\
$n=50$ & Biweight  &0.496 &0.542 &  &1.234  &1.194  & &1.228 & 1.228\\
& &(0.178) &(0.202) & &(0.190) &(0.206) & &(0.202) &(0.204)\\
& & & & & & & & &\\
 & Triweight  &0.522 &0.578  & &1.210  &1.192  & &1.230  &1.218 \\
& &(0.182) &(0.198) & &(0.216) &(0.216) & &(0.200) &(0.200)\\
& & & & & & & & &\\
& Triangular    	&0.524 &0.562  & &1.240  & 1.182  & &1.236  &  1.236\\
& &(0.178) &(0.192) & &(0.206) &(0.222) & &(0.210) &(0.204)\\
& & & & & & & & &\\
 & & & & & & & & &\\
\hline
& & & & & & & & &\\
 &Gaussian&  0.396&0.328& &1.266   & 1.182& &1.224& 1.216\\
& &(0.142) &(0.112) & &(0.178) &(0.174) & &(0.196) &(0.178)\\
& & & & & & & & &\\
 & Epanechnikov 	&0.346 &0.344 &  &1.240  &1.178  & & 1.222  &  1.214\\
& &(0.126) &(0.126) & &(0.174) &(0.176) & &(0.178) &(0.184)\\
& & & & & & & & &\\
$n=100$ & Biweight  &0.332&0.354 &  & 1.220  &1.170  & &1.224  & 1.218\\
& &(0.116) &(0.124) & &(0.182) &(0.190) & &(0.184) &(0.186)\\
& & & & & & & & &\\
 & Triweight  &0.346 &0.370  & &1.222  &1.180  & &1.236  &1.226 \\
& &(0.124) &(0.134) & &(0.186) &(0.182) & &(0.186) &(0.176)\\
& & & & & & & & &\\
& Triangular    	&0.340 &0.358  & &1.238  & 1.178  & &1.230  &  1.230\\
& &(0.118) &(0.122) & &(0.172) &(0.180) & &(0.186) &(0.176)\\
 & & & & & & & & &\\
\hline\hline
	\end{tabular}
\end{table}

\begin{table}[htbp]
\centering 
\caption{Average and standard deviation of $\mathcal{D}$ over $500$  replicates, with sample size $n=200,400$.}
	\label{tab2} 
\medskip
	\begin{tabular}{cccccccccc}
\hline\hline
& & & & & & & & &\\
Sample size & Kernel& \multicolumn{2}{c}{Model 1}& & \multicolumn{2}{c}{Model 2}& & \multicolumn{2}{c}{Model 3}\\
\cline{3-4}\cline{6-7}\cline{9-10}
& & $k$-NN& Kernel& &$k$-NN &Kernel & & $k$-NN&Kernel\\
\hline
& & & & & & & & &\\
 &Gaussian&  0.266&0.220& &1.146   & 1.174& &1.228& 1.224\\
& &(0.092) &(0.074) & &(0.186) &(0.138) & &(0.160) &(0.132)\\
& & & & & & & & &\\
 & Epanechnikov 	& 0.230 &0.226 &  &1.228  &1.170  & & 1.218  &  1.226\\
& &(0.083) &(0.080) & &(0.146) &(0.146) & &(0.158) &(0.142)\\
& & & & & & & & &\\
$n=200$ & Biweight  &0.238 &0.240 &  &1.222  &1.172  & &1.208 & 1.218\\
& &(0.082) &(0.084) & &(0.160) &(0.146) & &(0.158) &(0.150)\\
& & & & & & & & &\\
 & Triweight  &0.236 &0.244  & &1.214  &1.160  & &1.224  &1.226 \\
& &(0.080) &(0.086) & &(0.148) &(0.142) & &(0.156) &(0.154)\\
& & & & & & & & &\\
& Triangular    	&0.244&0.244  & &1.228  & 1.162  & &1.236  &  1.234\\
& &(0.085) &(0.082) & &(0.160) &(0.150) & &(0.142) &(0.134)\\
& & & & & & & & &\\
 & & & & & & & & &\\
\hline
& & & & & & & & &\\
 &Gaussian&  0.176&0.154& &1.266   & 1.184& &1.220& 1.220\\
& &(0.066) &(0.054) & &(0.122) &(0.096) & &(0.120) &(0.096)\\
& & & & & & & & &\\
 & Epanechnikov 	&0.164 &0.162 &  &1.228  &1.178  & & 1.218  &  1.222\\
& &(0.060) &(0.060) & &(0.116) &(0.106) & &(0.112) &(0.108)\\
& & & & & & & & &\\
$n=400$ & Biweight  &0.160&0.160 &  & 1.224  &1.174  & &1.220  & 1.222\\
& &(0.061) &(0.062) & &(0.116) &(0.106) & &(0.112) &(0.112)\\
& & & & & & & & &\\
 & Triweight  &0.162 &0.166  & &1.216  &1.168  & &1.224  &1.228 \\
& &(0.060) &(0.060) & &(0.114) &(0.104) & &(0.102) &(0.098)\\
& & & & & & & & &\\
& Triangular    	&0.160 &0.160  & &1.226  & 1.174  & &1.222  &  1.222\\
& &(0.058) &(0.058) & &(0.114) &(0.104) & &(0.102) &(0.098)\\
 & & & & & & & & &\\
\hline\hline
	\end{tabular}
\end{table}
\section{Proofs}\label{proofs}
\noindent We first give some preliminary lemmas. Then, these lemmas are used for proving Theorem \ref{theo1}. The proof of Theorem \ref{theo2} is identical to that of Theorem 2 in \cite{nkou22}; it is, therefore, omitted.
\subsection{Preliminary lemmas}

\bigskip

\noindent For $(\ell , j)\in\{1,\ldots,d\}^2$, we consider:
\begin{align}\label{u1}
 U_{n,\ell,j}^{(1)}& =\frac{1}{\sqrt{n}}\sum_{i=1}^n\Big\{g_\ell(Y_i)\left(\widehat{g}_{j , n}(Y_i)-g_{j}(Y_i)\right) + g_{j}(Y_i)\left(\widehat{g}_{\ell , n}(Y_i)-g_\ell(Y_i)\right)\Big\} \frac{\widehat{f}_{b_n}^2(Y_i)-f_{b_n}^2(Y_i)}{\widehat{f}_{b_n}^2(Y_i) f_{b_n}^2(Y_i)},
\end{align}
\begin{align}\label{u2}
U_{n,\ell,j}^{(2)}& = \frac{1}{\sqrt{n}}\sum_{i=1}^n \frac{\left(\widehat{g}_{\ell , n}(Y_i)-g_\ell(Y_i)\right)\left(\widehat{g}_{j , n}(Y_i)-g_{j}(Y_i)\right)}{\widehat{f}_{b_n}^2(Y_i)},
\end{align}
\begin{align}\label{u3}
U_{n,\ell,j}^{(3)} &= \frac{1}{\sqrt{n}}\sum_{i=1}^n R_{b_n,\ell}(Y_i)R_{b_n,j}(Y_i)\frac{\left(\widehat{f}_{b_n}^2(Y_i)-f_{b_n}^2(Y_i)\right)^2}{\widehat{f}_{b_n}^2(Y_i)f_{b_n}^2(Y_i)},
\end{align}
\begin{align}\label{u4}
U_{n,\ell,j}^{(4)} &= \frac{1}{\sqrt{n}}\sum_{i=1}^n\left(\widehat{f}_{b_n}(Y_i)-f_{b_n}(Y_i)\right)^2\frac{R_{b_n,\ell}(Y_i)R_{b_n,j}(Y_i)}{f_{b_n}^2(Y_i)},
\end{align}
where $R_{b_n,\ell}$ is the function defined in \eqref{rbnl}. Then, we have:

\bigskip

\begin{Lem}\label{cl14}
Under the assumptions \ref{as7} to \ref{as3}, \ref{as4}, \ref{as5}  and \ref{as6}, we have for  $(\ell,j)\in\{1,\ldots,d\}^2$ and  $m\in\{1,\ldots,4\}$,
$\left|U^{(m)}_{n,\ell,j}\right| = \mathrm{o}_{p}(1)$.
\end{Lem}
\noindent\textit{Proof.}
Arguing as in the proof of Lemma 5.2 in \cite{nkou22} and using Theorem 1 and Theorem 2 of \cite{bengono24}, we get, almost surely, 
\begin{eqnarray*}
\left|U_{n,\ell,j}^{(1)}\right|
                      & \leqslant &C_1C_2\, \rho_n\tau_n\,\sqrt{n}\,b_n^{-2}\,\left(C_2b_n^{-1}\tau_n + 2\, \right) \frac{1}{n}\sum_{i=1}^n\bigg(\left|R_\ell(Y_i)\right|+\left|R_{j}(Y_i)\right|\bigg),
\end{eqnarray*}
where $C_1$ and  $C_2$ are positive constants,   $\rho_n=\frac{k_{n}^4}{n^4}+\frac{M_n\sqrt{n\log(n)}}{k_{n}}$  and  $\tau_n=\frac{k_{n}^4}{n^4}+\frac{\sqrt{n\log(n)}}{k_{n}}$. Since, from Assumptions \ref{as7},  \ref{as5} and \ref{as6},
\begin{equation}\label{equivrt}
\rho_n\sim \frac{M_nn^{1/2}\log^{1/2}(n)}{k_n}\sim n^{1/2-c_1}\log(n),\,\,\,
\tau_n\sim \frac{n^{1/2}\log^{1/2}(n)}{k_n}\sim n^{1/2-c_1}\log^{1/2}(n),
\end{equation}
it follows that
\[
b_n^{-1}\tau_n \sim  n^{c_2+1/2-c_1}\log^{1/2}(n),\,\,\,\rho_n \tau_n\,\sqrt{n}\,b_n^{-2}\sim n^{2\left(3/4+c_2-c_1\right)}\log^{3/2}(n)
\]
and, consequently, that  $\lim_{n\rightarrow \infty}\rho_n \tau_n\,\sqrt{n}\,b_n^{-2}\,\left(b_n^{-1}\tau_n + 2\, \right) =0$ because $c_1-c_2>3/4>1/2$. Then, using   the law of large numbers  we deduce that  $
\left|U_{n,\ell,j}^{(1)}\right| = o_p(1)$. Similarly, as in the proof of Lemma 5.2 in \cite{nkou22}, we have, almost surely, the inequalities
\[
\left|U_{n,\ell,j}^{(2)}\right|\leqslant C_1^2n^{1/2}\,b_n^{-2}\,\rho_n^2,
\]
\[
\left|U_{n,\ell,j}^{(3)}\right| \leqslant C_2^2 \tau_n^2\,\sqrt{n}\,b_n^{-2}\,\left(C_2^2b_n^{-2}\tau_n^2 +4C_2b_n^{-1}\tau_n+ 4\right) \frac{1}{n}\sum_{i=1}^n\left|R_\ell(Y_i)R_{j}(Y_i)\right|,
\]
\[
\left|U_{n,\ell,j}^{(4)}\right|   \leqslant C_2^2\left(n^{1/2}b_n^{-2}\,\tau_n^2 \right) \frac{1}{n}\sum_{i=1}^n \left|R_\ell(Y_i)R_{j}(Y_i)\right|, 
\]
from which we deduce that $\left|U_{n,\ell,j}^{(2)}\right|=\mathrm{o}_{p}(1)$, $\left|U_{n,\ell,j}^{(3)}\right|=\mathrm{o}_{p}(1)$ and $\left|U_{n,\ell,j}^{(4)}\right|=\mathrm{o}_{p}(1)$ since, as above, we have  $\lim_{n\rightarrow \infty}\left( n^{1/2}\,b_n^{-2}\,\rho_n^2\right)=0$, $\lim_{n\rightarrow \infty}\left( \tau_n^2\,\sqrt{n}\,b_n^{-2}\,\left(C_2^2b_n^{-2}\tau_n^2 +4C_2b_n^{-1}\tau_n+ 4\right) \right)=0$  and $\lim_{n\rightarrow \infty}\left( n^{1/2}b_n^{-2}\,\tau_n^2  \right)=0$.
 \hfill $\Box$

\bigskip

\noindent Considering a sequence $(\beta_n)_{n\in\mathbb{N}^\ast}$ in $]0,1[$ such that $1-\beta_n\sim n^{-4}$, we put 
\begin{equation}\label{dnmp}
D_{n}^{-}(y)=\frac{k_{n}\sqrt{\beta_{n}}}{nf(y)}, \,\,\,D_{n}^{+}(y)=\frac{k_{n}}{n\sqrt{\beta_{n}}f(y)},
\end{equation}
and we consider, for $\ell\in\{1,\ldots,d\}$, the functions defined as
\begin{equation}\label{g1m}
\widehat{g}_{1,\ell,n}^{-}(y)=\frac{1}{nD_{n}^{+}(y)}\sum\limits_{i=1}^{n}X_{i\ell}\mathbf{1}_{\left\{X_{i\ell}\geqslant 0\right\}}K\left(\frac{Y_{i}-y}{D_{n}^{-}(y)}\right)
\end{equation}
and
\begin{equation}\label{g1p}
\widehat{g}_{1,\ell,n}^{+}(y)=\frac{1}{nD_{n}^{-}(y)}\sum\limits_{i=1}^{n}X_{i\ell}\mathbf{1}_{\left\{X_{i\ell}\geqslant 0\right\}}K\left(\frac{Y_{i}-y}{D_{n}^{+}(y)}\right).
\end{equation}
Then, we have:

\bigskip

\begin{Lem}\label{l3}
Under the assumptions \ref{as3}  to   \ref{as6}, we have for  $(\ell,j)\in\{1,\ldots,d\}^2$:
\begin{equation*}\label{resl3}
 \sqrt{n}\bigg(\mathbb{E}\bigg[\frac{g_{1,\ell}(Y)}{f_{b_n}^2(Y)}\Big(\widehat{g}^{-}_{1,j,n}(Y)-\widehat{g}^{+}_{1,j,n}(Y)\Big)\bigg]\bigg)=\mathrm{o}(1),
\end{equation*}
where $g_{1,\ell}$ is a function defined in \eqref{g1lg2l}.
\end{Lem}
\noindent\textit{Proof.} We have:
\begin{align*}
&\sqrt{n}\mathbb{E}\bigg[\frac{g_{1,\ell}(Y)}{f_{b_n}^2(Y)}\widehat{g}^{-}_{1,j,n}(Y)\bigg]\\
&=\frac{n^{3/2}\sqrt{\beta_{n}}}{k_{n}}\mathbb{E}\left[\frac{1}{n}\sum\limits_{i=1}^{n}\frac{g_{1,\ell}(Y)f(Y)}{f_{b_n}^2(Y)}X_{ij}\mathbf{1}_{\left\{X_{ij}\geqslant 0\right\}}K\left(\frac{nf(Y)(Y_{i}-Y)}{k_{n}\sqrt{\beta_{n}}} \right)\right]\notag\\
&=\frac{n^{3/2}\sqrt{\beta_{n}}}{k_{n}}\mathbb{E}\left[\frac{g_{1,\ell}(Y)f(Y)}{f_{b_n}^2(Y)}X_{1j}\mathbf{1}_{\left\{X_{1j}\geqslant 0\right\}}K\left(\frac{nf(Y)(Y_{1}-Y)}{k_{n}\sqrt{\beta_{n}}} \right)\right]\notag\\
&=\frac{n^{3/2}\sqrt{\beta_{n}}}{k_{n}}\int\int\int \frac{g_{1,\ell}(y)f(y)}{f_{b_n}^2(y)}x\mathbf{1}_{\mathbb{R}_+}(x)K\left(\frac{nf(y)(z-y)}{k_{n}\sqrt{\beta_{n}}} \right)\,f_{(X_{1j},Y_1,Y)}(x,z,y)\,\,dx\,dz\,dy\notag\\
&=\frac{n^{3/2}\sqrt{\beta_{n}}}{k_{n}}\int\int\int \frac{g_{1,\ell}(y)f(y)}{f_{b_n}^2(y)}x\mathbf{1}_{\mathbb{R}_+}(x)K\left(\frac{nf(y)(z-y)}{k_{n}\sqrt{\beta_{n}}} \right)\,f_{(X_{1j},Y_1)}(x,z)f(y)\,\,dx\,dz\,dy\notag\\
&=\frac{n^{3/2}\sqrt{\beta_{n}}}{k_{n}}\int\int\int \frac{g_{1,\ell}(y)f(y)}{f_{b_n}^2(y)}x\mathbf{1}_{\mathbb{R}_+}(x)K\left(\frac{nf(y)(z-y)}{k_{n}\sqrt{\beta_{n}}} \right)\,f_{(X_{j},Y)}(x,z)f(y)\,\,dx\,dz\,dy\notag\\
&=\frac{n^{3/2}\sqrt{\beta_{n}}}{k_{n}}\int\int\frac{g_{1,\ell}(y)f(y)}{f_{b_n}^2(y)}\bigg(\int x\mathbf{1}_{\mathbb{R}_+}(x)f_{(X_j,Y)}(x,z)\,dx\bigg)K\left(\frac{nf(y)(z-y)}{k_{n}\sqrt{\beta_{n}}} \right)f(y)\,\,dz\,dy\notag\\
&=\frac{n^{3/2}\sqrt{\beta_{n}}}{k_{n}}\int\int\frac{g_{1,\ell}(y)f(y)}{f_{b_n}^2(y)}g_{1,j}(z) K\left(\frac{nf(y)(z-y)}{k_{n}\sqrt{\beta_{n}}} \right)f(y)\,\,dz\,dy\notag.
\end{align*}
Since
\begin{eqnarray*}
\int g_{1,j}(z) K\left(\frac{nf(y)(z-y)}{k_{n}\sqrt{\beta_{n}}} \right) \,\,dz&=&\frac{k_n\sqrt{\beta_n}}{nf(y)}\int g_{1,j}\left(y+\frac{k_{n}\sqrt{\beta_{n}}}{nf(y)}t\right) K\left(t \right) \,\,dt,
\end{eqnarray*}
it follows
\begin{align*}
\sqrt{n}\mathbb{E}\bigg[\frac{g_{1,\ell}(Y)}{f_{b_n}^2(Y)}\widehat{g}^{-}_{1,j,n}(Y)\bigg]&=\sqrt{n} \beta_{n}\int\int\frac{g_{1,\ell}(y)}{f_{b_n}^2(y)}g_{1,j}\left(y+\frac{k_{n}\sqrt{\beta_{n}}}{nf(y)}t\right)  K\left(t \right)f(y) \,\,dt\,dy\notag\\
&=A_n+B_n,
\end{align*}
where
\[
A_n=\sqrt{n} \beta_{n}\int\int\frac{g_{1,\ell}(y)}{f_{b_n}^2(y)}\bigg(g_{1,j}\left(y+\frac{k_{n}\sqrt{\beta_{n}}}{nf(y)}t\right)-g_{1,j}\left(y\right)\bigg)  K\left(t \right)f(y) \,\,dt\,dy
\]
and
\[
B_n=\sqrt{n} \beta_{n}\int\int\frac{g_{1,\ell}(y)g_{1,j}(y)}{f_{b_n}^2(y)}    K\left(t \right)f(y) \,\,dt\,dy=\sqrt{n} \beta_{n}\int\frac{g_{1,\ell}(y)g_{1,j}(y)}{f_{b_n}^2(y)}    f(y) \,\,dy.
\]
By Taylor expansion, we have: 
\begin{equation*}
g_{1,j}\left(y+\frac{k_{n}\sqrt{\beta_{n}}}{nf(y)}t\right)=g_{1,j}(y)+\sum\limits_{m=1}^{2}\frac{g^{(m)}_{1,j}(y)}{m!} \left(\frac{k_{n}\sqrt{\beta_{n}}}{nf(y)}\right)^{m}t^{m}+\frac{1}{6}\left(\frac{k_{n}\sqrt{\beta_{n}}}{nf(y)}\right)^{3}t^{3}g^{(3)}_{1,j}\left(y+\theta\frac{k_{n}\sqrt{\beta_{n}}}{nf(y)}t\right),
\end{equation*}
where $0<\theta<1$. Thus 
\begin{eqnarray*}
\int \bigg(g_{1,j}\left(y+\frac{k_{n}\sqrt{\beta_{n}}}{nf(y)}t\right)-g_{1,j}\left(y\right)\bigg)  K\left(t \right) \,\,dt
&=&\sum\limits_{m=1}^{2}\frac{g^{(m)}_{1,j}(y)}{m!} \left(\frac{k_{n}\sqrt{\beta_{n}}}{nf(y)}\right)^{m}\int t^{m} K\left(t \right) \,\,dt \\
&&+\frac{1}{6}\left(\frac{k_{n}\sqrt{\beta_{n}}}{nf(y)}\right)^{3}\int t^{3}g^{(3)}_{1,j}\left(y+\theta\frac{k_{n}\sqrt{\beta_{n}}}{nf(y)}t\right) K(t)\,dt,
\end{eqnarray*}
and since $K$ is of order 3 (see Assumption \ref{as4}-$(iv)$), it follows
\begin{align*}
&\int \bigg(g_{1,j}\left(y+\frac{k_{n}\sqrt{\beta_{n}}}{nf(y)}t\right)-g_{1,j}\left(y\right)\bigg)  K\left(t \right) \,\,dt\\
&=\frac{1}{6}\left(\frac{k_{n}\sqrt{\beta_{n}}}{nf(y)}\right)^{3}\int t^{3}g^{(3)}_{1,j}\left(y+\theta\frac{k_{n}\sqrt{\beta_{n}}}{nf(y)}t\right) K(t)\,dt\\
&=\frac{1}{6}\left(\frac{k_{n}\sqrt{\beta_{n}}}{nf(y)}\right)^{3}\int t^{3}\bigg(g^{(3)}_{1,j}\left(y+\theta\frac{k_{n}\sqrt{\beta_{n}}}{nf(y)}t\right)-g^{(3)}_{1,j}\left(y\right)\bigg) K(t)\,dt.
\end{align*}
Hence
\[
A_n= \frac{k_{n}^3 \beta_{n}^{5/2}}{6n^{5/2}}\int\int\frac{g_{1,\ell}(y)t^3}{f_{b_n}^2(y)f^2(y)}\bigg(g^{(3)}_{1,j}\left(y+\theta\frac{k_{n}\sqrt{\beta_{n}}}{nf(y)}t\right)-g^{(3)}_{1,j}\left(y\right)\bigg)  K\left(t \right) \,\,dt\,dy
\]
and, using Assumption  \ref{as3},
\begin{eqnarray*}
\left\vert A_{n}\right\vert&\leqslant&c\frac{k_{n}^4 \beta_{n}^{3}}{6n^{7/2}}\theta\int t^4\vert K(t)\vert\,dt\times\int\frac{\vert g_{1,\ell}(y)\vert }{f_{b_n}^2(y)f^3(y)}\,dy\\
&\leqslant&c\frac{k_{n}^4 \beta_{n}^{3}}{6n^{7/2}}\int t^4\vert K(t)\vert\,dt\times\int\frac{\vert g_{1,\ell}(y)\vert }{f_{b_n}^2(y)f^3(y)}\,dy\\
&\leqslant&c\frac{k_{n}^4 \beta_{n}^{3}}{6n^{7/2}}\int t^4\vert K(t)\vert\,dt\times\int\frac{\vert R_{1,\ell}(y)\vert }{f_{b_n}(y)f^3(y)}\,dy\\
&\leqslant&c\frac{k_{n}^4 b_{n}^{-1}}{6c_0^4n^{7/2}}\int t^4\vert K(t)\vert\,dt\times\int\vert R_{1,\ell}(y)\vert f(y)\,dy\\
&=&c\frac{k_{n}^4 b_{n}^{-1}}{6c_0^4n^{7/2}}\int t^4\vert K(t)\vert\,dt\times \mathbb{E}\left( \vert R_{1,\ell}(Y)\vert\right). 
\end{eqnarray*}
From Assumptions  \ref{as5} and  \ref{as6},  $n^{-7/2}b_{n}^{-1}k_{n}^{4}\sim n^{4c_{1}+c_{2}-7/2}\rightarrow 0$ as  $n\rightarrow \infty$. Then, we deduce from the above inequality that $A_n=\mathrm{o}(1)$. Moreover, since
\begin{eqnarray*}
B_{n}
&= &\sqrt{n} \beta_{n}\int\frac{g_{1,\ell}(y)g_{1,j}(y)}{f_{b_n}^2(y)}    f(y) \,\,dy\\
&= &\sqrt{n} \beta_{n}\int  R_{1,\ell}(y)R_{1,j}(y)\varepsilon_{b_{n}}^{2}(y)   f(y) \,\,dy\\
&=&\sqrt{n} \beta_{n}\mathbb{E}\bigg[ R_{1,\ell}(Y)R_{1,j}(Y)\varepsilon_{b_{n}}^{2}(Y) \bigg],
\end{eqnarray*}
where 
\begin{equation}\label{eps}
\varepsilon_{b_{n}}(y)=\frac{f(y)}{f_{b_{n}}(y)},
\end{equation}
it follows 
\[
\left|B_n-\sqrt{n} \beta_{n}\mathbb{E}\left(R_{1,\ell}(Y)R_{1,j}(Y)\right)\right| \leqslant \sqrt{n} \beta_{n}\mathbb{E}\bigg[\vert R_{1,\ell}(Y)R_{1,j}(Y)\vert\,\,\vert 1-\varepsilon_{b_{n}}^{2}(Y)\vert \bigg].
\]
However, since $0\leqslant \varepsilon_{b_{n}}(y)\leqslant 1$, we have
\begin{equation*}
0\leqslant 1-\varepsilon_{b_{n}}^{2}(Y)=\left(1-\frac{f^{2}(Y)}{f^{2}_{b_{n}}(Y)}\right)\mathbf{1}_{\left\{f(Y)<b_{n}\right\}}\leqslant \mathbf{1}_{\left\{f(Y)<b_{n}\right\}}
\end{equation*}
and, therefore,
\[
\left|B_n-\sqrt{n} \beta_{n}\mathbb{E}\left(R_{1,\ell}(Y)R_{1,j}(Y)\right)\right| \leqslant \sqrt{n} \mathbb{E}\left[\vert R_{1,\ell}(Y)R_{1,j}(Y)\vert\,\,\mathbf{1}_{\left\{f(Y)<b_{n}\right\}}\right]=\mathrm{o}(1),
\] 
which implies that
\begin{equation*}
B_{n}=\sqrt{n} \beta_{n}\mathbb{E}\left(R_{1,\ell}(Y)R_{1,j}(Y)\right)+\mathrm{o}(1).
\end{equation*}
and, consequently, that
\begin{equation*}
\sqrt{n}\mathbb{E}\bigg[\frac{g_{1,\ell}(Y)}{f_{b_n}^2(Y)}\widehat{g}^{-}_{1,j,n}(Y)\bigg]=\sqrt{n} \beta_{n}\mathbb{E}\left(R_{1,\ell}(Y)R_{1,j}(Y)\right)+\mathrm{o}(1).
\end{equation*}
From this latter equality, we get
\begin{equation*}
\sqrt{n}\mathbb{E}\bigg[\frac{g_{1,\ell}(Y)}{f_{b_n}^2(Y)}\widehat{g}^{-}_{1,j,n}(Y)\bigg]-\sqrt{n} \mathbb{E}\left(R_{1,\ell}(Y)R_{1,j}(Y)\right)=-\sqrt{n} (1-\beta_{n})\mathbb{E}\left(R_{1,\ell}(Y)R_{1,j}(Y)\right)+\mathrm{o}(1),
\end{equation*}
and since $\sqrt{n}(1-\beta_n)\sim n^{-7/2}\rightarrow 0$ as $n\rightarrow\infty$, it follows that
\begin{equation}\label{res1}
\sqrt{n}\mathbb{E}\bigg[\frac{g_{1,\ell}(Y)}{f_{b_n}^2(Y)}\widehat{g}^{-}_{1,j,n}(Y)\bigg]=\sqrt{n} \mathbb{E}\left(R_{1,\ell}(Y)R_{1,j}(Y)\right)+\mathrm{o}(1).
\end{equation}
Furthermore, since $\widehat{g}^{+}_{1,j,n}$ is obtained from $\widehat{g}^{-}_{1,j,n}$ by replacing $\beta_n$ by $1/{\beta_n}$, we obtain  analogously,  
\begin{equation}\label{anag}
\sqrt{n}\mathbb{E}\bigg[\frac{g_{1,\ell}(Y)}{f_{b_n}^2(Y)}\widehat{g}^{+}_{1,j,n}(Y)\bigg]=A^\prime_n+B^\prime_n, 
\end{equation}
where
\begin{eqnarray}\label{anp}
\left\vert A_{n}^\prime\right\vert&\leqslant&c\frac{k_{n}^4}{6n^{7/2} \beta_{n}^{3}}\theta\int t^4\vert K(t)\vert\,dt\times\int\frac{\vert g_{1,\ell}(y)\vert }{f_{b_n}^2(y)f^3(y)}\,dy\nonumber\\
&\leqslant&c\frac{k_{n}^4 b_{n}^{-1}}{6c_0^4 \beta_{n}^{3}n^{7/2}}\int t^4\vert K(t)\vert\,dt\times \mathbb{E}\left( \vert R_{1,\ell}(Y)\vert\right). 
\end{eqnarray}
and
\begin{equation}\label{bnp}
\left|B_n^\prime-\frac{\sqrt{n}}{\beta_{n}}\mathbb{E}\left(R_{1,\ell}(Y)R_{1,j}(Y)\right)\right| \leqslant \frac{\sqrt{n}}{\beta_{n}} \mathbb{E}\left[\vert R_{1,\ell}(Y)R_{1,j}(Y)\vert\,\,\mathbf{1}_{\left\{f(Y)<b_{n}\right\}}\right].
\end{equation}
Since $\beta_n\rightarrow 1$ as $n\rightarrow\infty$, it follows from \eqref{anp},  \eqref{bnp} and Assumption  \ref{as9} that $A_n^\prime=\mathrm{o}(1)$,    $\left|B_n^\prime-\frac{\sqrt{n}}{\beta_{n}}\mathbb{E}\left(R_{1,\ell}(Y)R_{1,j}(Y)\right)\right|=\mathrm{o}(1)$ and, consequently, that    $B_n^\prime=\frac{\sqrt{n}}{\beta_{n}}\mathbb{E}\left(R_{1,\ell}(Y)R_{1,j}(Y)\right)+\mathrm{o}(1)$.
This yields
\[
B_n^\prime- \sqrt{n} \mathbb{E}\left(R_{1,\ell}(Y)R_{1,j}(Y)\right)=\frac{\sqrt{n}(1-\beta_n)}{\beta_n}\mathbb{E}\left(R_{1,\ell}(Y)R_{1,j}(Y)\right)+\mathrm{o}(1),
\]
and since
\[
\frac{\sqrt{n}(1-\beta_n)}{\beta_n}\sim \sqrt{n}\,n^{-4}=n^{-7/2},
\]
it follows that $B_n^\prime- \sqrt{n} \mathbb{E}\left(R_{1,\ell}(Y)R_{1,j}(Y)\right)=\mathrm{o}(1)$ and, from \eqref{anag}, that
\begin{equation}\label{res2}
\sqrt{n}\mathbb{E}\bigg[\frac{g_{1,\ell}(Y)}{f_{b_n}^2(Y)}\widehat{g}^{+}_{1,j,n}(Y)\bigg]= \sqrt{n} \mathbb{E}\left(R_{1,\ell}(Y)R_{1,j}(Y)\right)+\mathrm{o}(1).
\end{equation}
Then, \eqref{resl3} is obtained from \eqref{res1} and \eqref{res2}.
\hfill $\Box$

\bigskip

\noindent Now,  we put
\begin{equation*}
\mathscr{V}_{i,\ell,j,n}^{-} = \frac{n\sqrt{\beta_n}f(Y_i)}{2k_n} \int\int \left(\frac{g_{1,\ell}(Y_i)}{f_{b_n}^2(Y_i)}x \mathbf{1}_{\mathbb{R}_+}(x) + \frac{g_{1,\ell}(z)}{f_{b_n}^2(z)}X_{ij} \mathbf{1}_{\{ X_{ij} \geqslant 0 \}}\right)K\bigg(\frac{nf(Y_i)(z-Y_i)}{k_n\sqrt{\beta_n}}\bigg)f_{(X_j,Y)}(x, z) \, dx \, dz,
\end{equation*}

\begin{equation*}
\mathscr{V}_{i,\ell,j,n}^{+} = \frac{nf(Y_i)}{2k_n\sqrt{\beta_n}} \int\int \left(\frac{g_{1,\ell}(Y_i)}{f_{b_n}^2(Y_i)}x \mathbf{1}_{\mathbb{R}_+}(x) + \frac{g_{1,\ell}(z)}{f_{b_n}^2(z)}X_{ij} \mathbf{1}_{\{ X_{ij} \geqslant 0 \}}\right)K\bigg(\frac{n\sqrt{\beta_n}f(Y_i)(z-Y_i)}{k_n}\bigg)f_{(X_j,Y)}(x, z) \, dx \, dz.
\end{equation*}

\bigskip

\noindent The following lemma is obtained by applying  similar argument than in several works in order to approximate a sum by a U-statistic (e.g., \cite{powell89}). Because of the length of its  proof we omit it.

\bigskip

\begin{Lem}\label{l45}
Under the assumptions \ref{as7}, \ref{as2}, \ref{as8},  \ref{as5} and \ref{as6} we have for  $(\ell,j)\in\{1,\ldots,d\}^2$:
\begin{equation*}
\begin{split}
(i)\,\,&\frac{1}{\sqrt{n}}\sum_{i=1}^n \left\{\frac{g_{1,\ell}(Y_i)}{f_{b_n}^2(Y_i)}\widehat{g}^{-}_{1,j,n}(Y_{i}) - \mathbb{E}\left[ \frac{g_{1,\ell}(Y)}{f_{b_n}^2(Y)}\widehat{g}^{-}_{1,j,n}(Y) \right] \right\}=\frac{1}{\sqrt{n}} \sum_{i=1}^n \bigg\{ \mathscr{V}_{i,\ell,j,n}^{-}  - \mathbb{E}\left( \mathscr{V}_{i,\ell,j,n}^{-} \right) \bigg\} + \mathrm{o}_{p}(1);
\end{split}
\end{equation*}
\begin{equation*}
\begin{split}
(ii)\,\,&\frac{1}{\sqrt{n}}\sum_{i=1}^n \left\{\frac{g_{1,\ell}(Y_i)}{f_{b_n}^2(Y_i)}\widehat{g}^{+}_{1,j,n}(Y_{i}) - \mathbb{E}\left[ \frac{g_{1,\ell}(Y)}{f_{b_n}^2(Y)}\widehat{g}^{+}_{1,j,n}(Y) \right] \right\}=\frac{1}{\sqrt{n}} \sum_{i=1}^n \bigg\{ \mathscr{V}_{i,\ell,j,n}^{+}  - \mathbb{E}\left( \mathscr{V}_{i,\ell,j,n}^{+} \right) \bigg\} + \mathrm{o}_{p}(1).
\end{split}
\end{equation*}
\end{Lem}

\bigskip

\noindent Then, using this lemma we obtain the following result.

\bigskip

\begin{Lem}\label{l67}
Under the assumptions  \ref{as0},  \ref{as8} and  \ref{as4} we have  for $(\ell,j)\in\{1,\ldots,d\}^2$:
\begin{align*}
(i)\,\,&\frac{1}{\sqrt{n}}\sum\limits_{i=1}^{n}\left\{\frac{g_{1,\ell}(Y_i)}{f_{b_n}^2(Y_i)}\widehat{g}^{-}_{1,j,n}(Y_{i})+\frac{g_{1,j}(Y_i)}{f_{b_n}^2(Y_i)}\widehat{g}^{-}_{1,\ell,n}(Y_{i})-\mathbb{E}\left[\frac{g_{1,\ell}(Y)}{f_{b_n}^2(Y)}\widehat{g}^{-}_{1,j,n}(Y)+\frac{g_{1,j}(Y)}{f_{b_n}^2(Y)}\widehat{g}^{-}_{1,\ell,n}(Y)\right]\right\}\notag\\
&=\frac{1}{\sqrt{n}}\sum\limits_{i=1}^{n}\bigg\{\frac{g_{1,\ell}(Y_i) g_{1,j}(Y_i)}{f_{b_n}^2(Y_i)}+\frac{g_{1,\ell}(Y_i)f(Y_i)}{2f_{b_n}^2(Y_i)}X_{ij}\mathbf{1}_{\left\{X_{ij}\geqslant 0\right\}}+\frac{g_{1,j}(Y_i)f(Y_i)}{2f_{b_n}^2(Y_i)}X_{i\ell}\mathbf{1}_{\left\{X_{i\ell}\geqslant 0\right\}}\notag\\
&\hspace{0.5cm}-\mathbb{E}\left[\frac{g_{1,\ell}(Y) g_{1,j}(Y)}{f_{b_n}^2(Y)}+\frac{g_{1,\ell}(Y)f(Y)}{2f_{b_n}^2(Y)}X_{ij}\mathbf{1}_{\left\{X_{ij}\geqslant 0\right\}}+\frac{g_{1,j}(Y)f(Y)}{2f_{b_n}^2(Y)}X_{i\ell}\mathbf{1}_{\left\{X_{i\ell}\geqslant 0\right\}}\right]\bigg{\}}+\mathrm{o}_{p}(1);
\end{align*}
\begin{align*}
(ii)\,\,&\frac{1}{\sqrt{n}}\sum\limits_{i=1}^{n}\left\{\frac{g_{1,\ell}(Y_i)}{f_{b_n}^2(Y_i)}\widehat{g}^{+}_{1,j,n}(Y_{i})+\frac{g_{1,j}(Y_i)}{f_{b_n}^2(Y_i)}\widehat{g}^{+}_{1,\ell,n}(Y_{i})-\mathbb{E}\left[\frac{g_{1,\ell}(Y)}{f_{b_n}^2(Y)}\widehat{g}^{+}_{1,j,n}(Y)+\frac{g_{1,j}(Y)}{f_{b_n}^2(Y)}\widehat{g}^{+}_{1,\ell,n}(Y)\right]\right\}\notag\\
&=\frac{1}{\sqrt{n}}\sum\limits_{i=1}^{n}\bigg\{\frac{g_{1,\ell}(Y_i) g_{1,j}(Y_i)}{f_{b_n}^2(Y_i)}+\frac{g_{1,\ell}(Y_i)f(Y_i)}{2f_{b_n}^2(Y_i)}X_{ij}\mathbf{1}_{\left\{X_{ij}\geqslant 0\right\}}+\frac{g_{1,j}(Y_i)f(Y_i)}{2f_{b_n}^2(Y_i)}X_{i\ell}\mathbf{1}_{\left\{X_{i\ell}\geqslant 0\right\}}\notag\\
&\hspace{0.5cm}-\mathbb{E}\left[\frac{g_{1,\ell}(Y) g_{1,j}(Y)}{f_{b_n}^2(Y)}+\frac{g_{1,\ell}(Y)f(Y)}{2f_{b_n}^2(Y)}X_{ij}\mathbf{1}_{\left\{X_{ij}\geqslant 0\right\}}+\frac{g_{1,j}(Y)f(Y)}{2f_{b_n}^2(Y)}X_{i\ell}\mathbf{1}_{\left\{X_{i\ell}\geqslant 0\right\}}\right]\bigg{\}}+\mathrm{o}_{p}(1);
\end{align*}
\end{Lem}
\noindent\textit{Proof.} 
Clearly,
\begin{eqnarray*}
\mathscr{V}_{i,\ell,j,n}^{-} &= &\frac{n\sqrt{\beta_n}g_{1,\ell}(Y_i)f(Y_i)}{2k_nf_{b_n}^2(Y_i)} \int\left(\int x \mathbf{1}_{\mathbb{R}_+}(x)\, f_{(X_j,Y)}(x, z) \, dx\right)K\bigg(\frac{nf(Y_i)(z-Y_i)}{k_n\sqrt{\beta_n}}\bigg) \, dz\\
& &+\frac{n\sqrt{\beta_n}f(Y_i)}{2k_n}X_{ij} \mathbf{1}_{\{ X_{ij} \geqslant 0 \}} \int \frac{g_{1,\ell}(z)}{f_{b_n}^2(z)}\left(\int f_{(X_j,Y)}(x, z) \, dx\right)K\bigg(\frac{nf(Y_i)(z-Y_i)}{k_n\sqrt{\beta_n}}\bigg) \, dz\\
 &= &\frac{n\sqrt{\beta_n}g_{1,\ell}(Y_i)f(Y_i)}{2k_nf_{b_n}^2(Y_i)} \int g_{1,j} (z)K\bigg(\frac{nf(Y_i)(z-Y_i)}{k_n\sqrt{\beta_n}}\bigg) \, dz\\
& &+\frac{n\sqrt{\beta_n}f(Y_i)}{2k_n}X_{ij} \mathbf{1}_{\{ X_{ij} \geqslant 0 \}} \int \frac{g_{1,\ell}(z)}{f_{b_n}^2(z)}f(z)K\bigg(\frac{nf(Y_i)(z-Y_i)}{k_n\sqrt{\beta_n}}\bigg) \, dz\\
 &= &\frac{n\sqrt{\beta_n}g_{1,\ell}(Y_i)f(Y_i)}{2k_nf_{b_n}^2(Y_i)} \int R_{1,j} (z)f(z)K\bigg(\frac{nf(Y_i)(z-Y_i)}{k_n\sqrt{\beta_n}}\bigg) \, dz\\
& &+\frac{n\sqrt{\beta_n}f(Y_i)}{2k_n}X_{ij} \mathbf{1}_{\{ X_{ij} \geqslant 0 \}} \int \frac{g_{1,\ell}(z)}{f_{b_n}^2(z)}f(z)K\bigg(\frac{nf(Y_i)(z-Y_i)}{k_n\sqrt{\beta_n}}\bigg) \, dz.
\end{eqnarray*}
Since
\begin{eqnarray*}
\frac{n\sqrt{\beta_n}f(Y_i)}{k_n}\int K\bigg(\frac{nf(Y_i)(z-y)}{k_n\sqrt{\beta_n}}\bigg) \,dz=\beta_n\int K(t)\,dt=\beta_n,
\end{eqnarray*}
it follows
\begin{equation}\label{decvm}
\mathscr{V}_{i,\ell,j,n}^{-} =\mathscr{C}_{\ell,j,n}(Y_i)+\mathscr{D}_{\ell,n}(X_{ij},Y_i)+\mathscr{E}_{\ell,j,n}(Y_i),
\end{equation}
where
\[
\mathscr{C}_{\ell,j,n}(Y_i)=\frac{n\sqrt{\beta_n}g_{1,\ell}(Y_i)f(Y_i)}{2k_nf_{b_n}^2(Y_i)} \int \bigg(R_{1,j} (z)f(z) - R_{1,j} (Y_i)f(Y_i) \bigg)K\bigg(\frac{nf(Y_i)(z-Y_i)}{k_n\sqrt{\beta_n}}\bigg) \, dz,
\]
\[
\mathscr{D}_{\ell,n}(X_{ij},Y_i)=\frac{n\sqrt{\beta_n}f(Y_i)}{2k_n}X_{ij} \mathbf{1}_{\{ X_{ij} \geqslant 0 \}} \int \bigg(\frac{g_{1,\ell}(z)}{f_{b_n}^2(z)}f(z)- \frac{g_{1,\ell}(Y_i)}{f_{b_n}^2(Y_i)}f(Y_i)\bigg)K\bigg(\frac{nf(Y_i)(z-Y_i)}{k_n\sqrt{\beta_n}}\bigg) \, dz
\]
and
\begin{eqnarray*}
\mathscr{E}_{\ell,j,n}(Y_i)&=&\frac{g_{1,\ell}(Y_i)R_{1,j} (Y_i)f(Y_i)\beta_n}{2f_{b_n}^2(Y_i)}+\frac{g_{1,\ell}(Y_i)f(Y_i)\beta_n}{2f_{b_n}^2(Y_i)}X_{ij} \mathbf{1}_{\{ X_{ij} \geqslant 0 \}}\\
& =&
\frac{1}{2}R_{1,\ell}(Y_i)R_{1,j} (Y_i)\varepsilon_{b_n}^2(Y_i)\beta_n +\frac{1}{2}R_{1,\ell}(Y_i)\varepsilon_{b_n}^2(Y_i)\beta_nX_{ij} \mathbf{1}_{\{ X_{ij} \geqslant 0 \}},
\end{eqnarray*}
$\varepsilon_{b_n}$ being defined in \eqref{eps}.
We have
\begin{align*}
\mathbb{E}\left(\mathscr{C}_{\ell,j,n}^2(Y_i)\right)&=\mathbb{E}\left[ \frac{n^2 \beta_n g_{1,\ell}^2(Y_i)f^2(Y_i)}{4k_n^2f_{b_n}^4(Y_i)} \bigg(\int \bigg(R_{1,j} (z)f(z) - R_{1,j} (Y_i)f(Y_i) \bigg)K\bigg(\frac{nf(Y_i)(z-Y_i)}{k_n\sqrt{\beta_n}}\bigg) \, dz\bigg)^2 \right]\notag\\
&=\int\frac{n^2 \beta_n g_{1,\ell}^2(y)f^2(y)}{4k_n^2f_{b_n}^4(y)} \bigg(\int \bigg(g_{1,j} (z) - g_{1,j} (y)\bigg)K\bigg(\frac{nf(y)(z-y)}{k_n\sqrt{\beta_n}}\bigg) \, dz\bigg)^2f(y)\,\,dy\notag\\
&=\int\frac{ \beta_n ^2g_{1,\ell}^2(y)}{4f_{b_n}^4(y)} \bigg(\int \bigg(g_{1,j} \left(y+\frac{ k_n\sqrt{\beta_n}}{nf(y) }t\right) - g_{1,j} (y)\bigg)K(t) \, dt\bigg)^2f(y)\,\,dy\notag\\
&\leqslant \int\frac{ \beta_n ^2g_{1,\ell}^2(y)}{4f_{b_n}^4(y)} \bigg(\int \bigg\vert g_{1,j} \left(y+\frac{ k_n\sqrt{\beta_n}}{nf(y) }t\right) - g_{1,j} (y)\bigg\vert\,\,\vert K(t)\vert \, dt\bigg)^2f(y)\,\,dy\\
&\leqslant \int\frac{c^2k_n^2 \beta_n ^3R_{1,\ell}^2(y)}{4n^2f_{b_n}^4(y)}f(y)\,\,dy\,\, \bigg(\int \vert t\vert\,\,\vert K(t)\vert \, dt\bigg)^2\\
&\leqslant\frac{ c^2k_n^2 n^{-2}b_n^{-2}}{4c_0^2}\mathbb{E}\left( R_{1,\ell}^2(Y)\varepsilon_{b_n}^2(Y)\right)\,\, \bigg(\int \vert t\vert\,\,\vert K(t)\vert \, dt\bigg)^2\\
&\leqslant\frac{ c^2k_n^2 n^{-2}b_n^{-2}}{4c_0^2}\mathbb{E}\left( R_{1,\ell}^2(Y)\right)\,\, \bigg(\int \vert t\vert\,\,\vert K(t)\vert \, dt\bigg)^2.
\end{align*}
Since $k_n^2 n^{-2}b_n^{-2}\sim n^{2c_{1}+2c_{2}-2}=n^{-2(1-c_{1}-c_{2})}$ and  $c_{1}+c_{2}<1$, we deduce that  $k_n^2 n^{-2}b_n^{-2}\rightarrow 0$ as $n\rightarrow\infty$, and from the above inequality that    $\mathbb{E}\left(\mathscr{C}_{\ell,j,n}^2(Y_i)\right)=\mathrm{o}(1)$. Since $Var\left(n^{-1/2}\sum_{i=1}^n\mathscr{C}_{\ell,j,n}(Y_i)\right) =Var\left( \mathscr{C}_{\ell,j,n}(Y)\right)\leqslant \mathbb{E} \bigg(\mathscr{C}_{\ell,j,n}^2(Y)\bigg) $, it follows from Bienaym\'e-Chebyshev  inequality that 
\begin{equation}\label{var1}
\frac{1}{\sqrt{n}} \sum_{i=1}^n\left(\mathscr{C}_{\ell,j,n}(Y_i)-\mathbb{E}  \left(\mathscr{C}_{\ell,j,n}(Y)\right) \right)=\mathrm{o}_{p}(1).
\end{equation} 
Furthermore, 
\begin{align*}
&\mathbb{E}\left[\mathscr{D}^2_{\ell,n}(X_{j},Y)\right]\\
&=\frac{n^2\beta_n}{4k_n^2}\mathbb{E}\left[ X_{j}^2 \mathbf{1}_{\{ X_{j} \geqslant 0 \}} \bigg(\int \bigg(\frac{g_{1,\ell}(z)}{f_{b_n}^2(z)}f(z)- \frac{g_{1,\ell}(Y)}{f_{b_n}^2(Y)}f(Y)\bigg)f(Y)K\bigg(\frac{nf(Y)(z-Y)}{k_n\sqrt{\beta_n}}\bigg) \, dz\bigg)^2\right]\\
&=\frac{n^2\beta_n}{4k_n^2}\mathbb{E}\left[ X_{j}^2 \mathbf{1}_{\{ X_{j} \geqslant 0 \}} \bigg(\int \bigg(R_{1,\ell}(z)\varepsilon_{b_n}^2(z)- R_{1,\ell}(Y)\varepsilon_{b_n}^2(Y)\bigg)f(Y)K\bigg(\frac{nf(Y)(z-Y)}{k_n\sqrt{\beta_n}}\bigg) \, dz\bigg)^2\right]\\
&=\frac{n^2\beta_n}{4k_n^2}\mathbb{E}\left[ \mathbb{E}\bigg(X_{j}^2 \mathbf{1}_{\{ X_{j} \geqslant 0 \}} \bigg(\int \bigg(R_{1,\ell}(z)\varepsilon_{b_n}^2(z)- R_{1,\ell}(Y)\varepsilon_{b_n}^2(Y)\bigg)f(Y)K\bigg(\frac{nf(Y)(z-Y)}{k_n\sqrt{\beta_n}}\bigg) \, dz\bigg)^2\bigg\vert Y\bigg)\right]\\
&=\frac{1}{4}\mathbb{E}\bigg[ \mathbb{E}\left(X_{j}^2 \mathbf{1}_{\{ X_{j} \geqslant 0 \}}\vert Y\right) \varphi^2_{\ell,n}(Y)\bigg],
\end{align*}
where
\begin{eqnarray*}
 \varphi_{\ell,n}(Y)&=&\frac{n\sqrt{\beta_n}}{k_n}\int \bigg(R_{1,\ell}(z)\varepsilon_{b_n}^2(z)- R_{1,\ell}(Y)\varepsilon_{b_n}^2(Y)\bigg)f(Y)K\bigg(\frac{nf(Y)(z-Y)}{k_n\sqrt{\beta_n}}\bigg) \, dz\\
&=&\beta_n\int \bigg(R_{1,\ell}\left(Y+\frac{k_n\sqrt{\beta_n}}{nf(Y)}t\right)\varepsilon_{b_n}^2\left(Y+\frac{k_n\sqrt{\beta_n}}{nf(Y)}t\right)- R_{1,\ell}(Y)\varepsilon_{b_n}^2(Y)\bigg)K(t)\, dt
\end{eqnarray*}
From the continuity of $f$ and $g_{1,\ell}$, which implies that of $R_{1,\ell}$ and  $\varepsilon_{b_n}$, and since $n^{-1}k_n\sim n^{c_1-1}\rightarrow 0$ and $\beta_n\rightarrow 1$ as $n\rightarrow \infty$, we obtain
\[
\lim_{n\rightarrow \infty}\bigg(R_{1,\ell}\left(Y+\frac{k_n\sqrt{\beta_n}}{nf(Y)}t\right)\varepsilon_{b_n}^2\left(Y+\frac{k_n\sqrt{\beta_n}}{nf(Y)}t\right)- R_{1,\ell}(Y)\varepsilon_{b_n}^2(Y)\bigg)K(t)=0.
\]
 Further, since  $0\leqslant \varepsilon_{b_n}\leqslant 1$,  it follows 
\begin{align*}
& \bigg\vert\bigg(R_{1,\ell}\left(Y+\frac{k_n\sqrt{\beta_n}}{nf(Y)}t\right)\varepsilon_{b_n}^2\left(Y+\frac{k_n\sqrt{\beta_n}}{nf(Y)}t\right)- R_{1,\ell}(Y)\varepsilon_{b_n}^2(Y)\bigg)K(t)\bigg\vert\nonumber\\
&\leqslant \bigg\vert \bigg(R_{1,\ell}\left(Y+\frac{k_n\sqrt{\beta_n}}{nf(Y)}t\right)- R_{1,\ell}(Y)\bigg)\varepsilon_{b_n}^2\left(Y+\frac{k_n\sqrt{\beta_n}}{nf(Y)}t\right)\nonumber\\
& + R_{1,\ell}(Y) \bigg(\varepsilon_{b_n}^2\left(Y+\frac{k_n\sqrt{\beta_n}}{nf(Y)}t\right)-\varepsilon_{b_n}^2\left(Y\right)\bigg)\bigg\vert\times \left\vert K(t)\right\vert\nonumber\\
&\leqslant\bigg(  \bigg\vert R_{1,\ell}\left(Y+\frac{k_n\sqrt{\beta_n}}{nf(Y)}t\right)- R_{1,\ell}(Y)\bigg\vert+2\vert R_{1,\ell}(Y)  \vert\bigg) \left\vert K(t)\right\vert\nonumber\\
&\leqslant c\frac{k_n\sqrt{\beta_n}}{nf(Y)}\vert t\vert \, \left\vert K(t)\right\vert +2\vert  R_{1,\ell}(Y)  \vert\, \left\vert K(t)\right\vert.
\end{align*}
Since $n^{-1}k_n \rightarrow 0$   as $n\rightarrow \infty$, we have for $n$ large enough
\begin{align}\label{dom1}
& \bigg\vert\bigg(R_{1,\ell}\left(Y+\frac{k_n\sqrt{\beta_n}}{nf(Y)}t\right)\varepsilon_{b_n}^2\left(Y+\frac{k_n\sqrt{\beta_n}}{nf(Y)}t\right)- R_{1,\ell}(Y)\varepsilon_{b_n}^2(Y)\bigg)K(t)\bigg\vert\nonumber\\
&\leqslant \frac{c}{f(Y)}\vert t\vert \, \left\vert K(t)\right\vert +2\vert   R_{1,\ell}(Y)  \vert\, \left\vert K(t)\right\vert\nonumber\\
&\leqslant \frac{c}{c_0}\vert t\vert \, \left\vert K(t)\right\vert +2\vert   R_{1,\ell}(Y)  \vert\, \left\vert K(t)\right\vert,
\end{align}
with
\begin{align*}
&\int  \frac{c}{c_0}\vert t\vert \, \left\vert K(t)\right\vert +2\vert   R_{1,\ell}(Y)  \vert\, \left\vert K(t)\right\vert\,\,dt=\frac{c}{c_0}\int\vert t\vert \, \left\vert K(t)\right\vert  dt+2\vert   R_{1,\ell}(Y)  \vert\, \int\left\vert K(t)\right\vert\,dt<\infty.
\end{align*}
Then,  using the dominated convergence theorem we get: $\lim_{n\rightarrow \infty}\left(\varphi_{\ell,n}(Y)\right)=0$. Moreover, using \eqref{dom1}, we have
\begin{eqnarray*}
\mathbb{E}\left(X_{j}^2 \mathbf{1}_{\{ X_{j} \geqslant 0 \}}\vert Y\right) \varphi^2_{\ell,n}(Y) &\leqslant&\mathbb{E}\left(X_{j}^2 \mathbf{1}_{\{ X_{j} \geqslant 0 \}}\vert Y\right) \bigg(\int \frac{c}{c_0}\vert t\vert \, \left\vert K(t)\right\vert +2\vert   R_{1,\ell}(Y)  \vert\, \left\vert K(t)\right\vert \,\,dt\bigg)^2\\
&\leqslant&\frac{2c^2}{c_0^2} \bigg(\int \vert t\vert \, \left\vert K(t)\right\vert\,\,dt\bigg)^2\mathbb{E}\left(X_{j}^2 \vert Y\right)+4 \bigg(\int \left\vert K(t)\right\vert\,\,dt\bigg)^2\mathbb{E}\left(X_{j}^2 \vert Y\right)R^2_{1,\ell}(Y).
\end{eqnarray*}
Since
\begin{align*}
& \mathbb{E}\bigg[ \frac{2c^2}{c_0^2} \bigg(\int \vert t\vert \, \left\vert K(t)\right\vert\,\,dt\bigg)^2\mathbb{E}\left(X_{j}^2 \vert Y\right)+4 \bigg(\int \left\vert K(t)\right\vert\,\,dt\bigg)^2\mathbb{E}\left(X_{j}^2 \vert Y\right)R^2_{1,\ell}(Y)\bigg]\\
  &= \frac{2c^2}{c_0^2} \bigg(\int \vert t\vert \, \left\vert K(t)\right\vert\,\,dt\bigg)^2\mathbb{E}\left(X_{j}^2 \right)+4 \bigg(\int \left\vert K(t)\right\vert\,\,dt\bigg)^2\mathbb{E}\left(X_{j}^2 R^2_{1,\ell}(Y)\right)\\
 & \leqslant \frac{2c^2}{c_0^2} \bigg(\int \vert t\vert \, \left\vert K(t)\right\vert\,\,dt\bigg)^2\mathbb{E}\left(X_{j}^2 \right)+4 \bigg(\int \left\vert K(t)\right\vert\,\,dt\bigg)^2\mathbb{E}^{1/2}\left(X_{j}^4\right)\, \mathbb{E}^{1/2}\left(R^4_{1,\ell}(Y)\right)<\infty,
\end{align*}
we can use again  the dominated convergence theorem which yields
\[
\lim_{n\rightarrow \infty}\mathbb{E}\left[\mathscr{D}^2_{\ell,n}(X_{j},Y)\right]=\frac{1}{4}\mathbb{E}\bigg[ \mathbb{E}\left(X_{j}^2 \mathbf{1}_{\{ X_{j} \geqslant 0 \}}\vert Y\right) \lim_{n\rightarrow \infty}\varphi^2_{\ell,n}(Y)\bigg]=0.
\]
Since $Var\left(n^{-1/2}\sum_{i=1}^n\mathscr{D}_{\ell,n}(X_{ij},Y_i)\right) =Var\left( \mathscr{D}_{\ell,n}(X_{j},Y)\right)\leqslant \mathbb{E} \left(\mathscr{D}^2_{\ell,n}(X_{j},Y)\right) $, we conclude, by using  Bienaym\'e-Chebyshev  inequality,  that
\begin{equation}\label{var2}
\frac{1}{\sqrt{n}} \sum_{i=1}^n\left(\mathscr{D}_{\ell,n}(X_{ij},Y_i)-\mathbb{E}  \left(\mathscr{D}_{\ell,n}(X_{j},Y)\right) \right)=\mathrm{o}_{p}(1).
\end{equation} 
Then, from \eqref{decvm}, \eqref{var1} and \eqref{var2} we obtain
\begin{align*}
&\frac{1}{\sqrt{n}}\sum\limits_{i=1}^n\left\{\mathscr{V}_{i,\ell,j,n}^{-}-\mathbb{E}\left(\mathscr{V}_{i,\ell,j,n}^{-}\right)\right\}\nonumber\\
&=\frac{1}{\sqrt{n}}\sum\limits_{i=1}^n\bigg\{\frac{1}{2}R_{1,\ell}(Y_i)R_{1,j} (Y_i)\varepsilon_{b_n}^2(Y_i)\beta_n +\frac{1}{2}R_{1,\ell}(Y_i)\varepsilon_{b_n}^2(Y_i)\beta_nX_{ij} \mathbf{1}_{\{ X_{ij} \geqslant 0 \}}\nonumber\\
&\hspace{0.2cm}-\mathbb{E}\left(\frac{1}{2}R_{1,\ell}(Y)R_{1,j} (Y)\varepsilon_{b_n}^2(Y)\beta_n +\frac{1}{2}R_{1,\ell}(Y)\varepsilon_{b_n}^2(Y)\beta_nX_{j} \mathbf{1}_{\{ X_{j} \geqslant 0 \}}\right)\bigg\}+\mathrm{o}_{p}(1),
\end{align*}
and from Lemma \ref{l45}  it follows
\begin{align*}
&\frac{1}{\sqrt{n}}\sum_{i=1}^n \left\{\frac{g_{1,\ell}(Y_i)}{f_{b_n}^2(Y_i)}\widehat{g}^{-}_{1,j,n}(Y_{i}) - \mathbb{E}\left[ \frac{g_{1,\ell}(Y)}{f_{b_n}^2(Y)}\widehat{g}^{-}_{1,j,n}(Y) \right] \right\}\\
&=\frac{1}{\sqrt{n}}\sum\limits_{i=1}^n\bigg\{\frac{1}{2}R_{1,\ell}(Y_i)R_{1,j} (Y_i)\varepsilon_{b_n}^2(Y_i)\beta_n +\frac{1}{2}R_{1,\ell}(Y_i)\varepsilon_{b_n}^2(Y_i)\beta_nX_{ij} \mathbf{1}_{\{ X_{ij} \geqslant 0 \}}\nonumber\\
&\hspace{0.2cm}-\mathbb{E}\left(\frac{1}{2}R_{1,\ell}(Y)R_{1,j} (Y)\varepsilon_{b_n}^2(Y)\beta_n +\frac{1}{2}R_{1,\ell}(Y)\varepsilon_{b_n}^2(Y)\beta_nX_{j} \mathbf{1}_{\{ X_{j} \geqslant 0 \}}\right)\bigg\}+\mathrm{o}_{p}(1)\\
&=\frac{1}{\sqrt{n}}\sum\limits_{i=1}^n\bigg\{\frac{1}{2}R_{1,\ell}(Y_i)R_{1,j} (Y_i)\varepsilon_{b_n}^2(Y_i)+\frac{1}{2}R_{1,\ell}(Y_i)\varepsilon_{b_n}^2(Y_i)X_{ij} \mathbf{1}_{\{ X_{ij} \geqslant 0 \}}\nonumber\\
&\hspace{0.2cm}-\mathbb{E}\left(\frac{1}{2}R_{1,\ell}(Y)R_{1,j} (Y)\varepsilon_{b_n}^2(Y)+\frac{1}{2}R_{1,\ell}(Y)\varepsilon_{b_n}^2(Y)X_{j} \mathbf{1}_{\{ X_{j} \geqslant 0 \}}\right)\bigg\}+\delta_n+\mathrm{o}_{p}(1),
\end{align*}
where
\begin{align*}
\delta_n&=\frac{\beta_n-1}{\sqrt{n}}\sum\limits_{i=1}^n\bigg\{\frac{1}{2}R_{1,\ell}(Y_i)R_{1,j} (Y_i)\varepsilon_{b_n}^2(Y_i )+\frac{1}{2}R_{1,\ell}(Y_i)\varepsilon_{b_n}^2(Y_i)X_{ij} \mathbf{1}_{\{ X_{ij} \geqslant 0 \}}\nonumber\\
&\hspace{0.2cm}-\mathbb{E}\left(\frac{1}{2}R_{1,\ell}(Y)R_{1,j} (Y)\varepsilon_{b_n}^2(Y)+\frac{1}{2}R_{1,\ell}(Y)\varepsilon_{b_n}^2(Y)X_{j} \mathbf{1}_{\{ X_{j} \geqslant 0 \}}\right)\bigg\},
\end{align*}
so that
\begin{align*}
\vert\delta_n\vert&\leqslant\frac{(1-\beta_n)\sqrt{n}}{2}\bigg\{\frac{1}{n}\sum\limits_{i=1}^n\left\vert R_{1,\ell}(Y_i)R_{1,j} (Y_i)\right\vert+\frac{1}{n}\sum\limits_{i=1}^n\left\vert R_{1,\ell}(Y_i)X_{ij} \right\vert+\mathbb{E}\left(\left\vert R_{1,\ell}(Y)R_{1,j} (Y)\right\vert+\left\vert R_{1,\ell}(Y)X_{j} \right\vert\right)\bigg\}.
\end{align*}
Since $(1-\beta_n)\sqrt{n}\sim n^{-7/2}$, we deduce from the preceding inequality and the strong law of large numbers that $\delta_n=\mathrm{o}_{p}(1)$ and, consequently, that
\begin{align}\label{res1i}
&\frac{1}{\sqrt{n}}\sum_{i=1}^n \left\{\frac{g_{1,\ell}(Y_i)}{f_{b_n}^2(Y_i)}\widehat{g}^{-}_{1,j,n}(Y_{i}) - \mathbb{E}\left[ \frac{g_{1,\ell}(Y)}{f_{b_n}^2(Y)}\widehat{g}^{-}_{1,j,n}(Y) \right] \right\}\nonumber\\
&=\frac{1}{\sqrt{n}}\sum\limits_{i=1}^n\bigg\{\frac{1}{2}R_{1,\ell}(Y_i)R_{1,j} (Y_i)\varepsilon_{b_n}^2(Y_i)+\frac{1}{2}R_{1,\ell}(Y_i)\varepsilon_{b_n}^2(Y_i)X_{ij} \mathbf{1}_{\{ X_{ij} \geqslant 0 \}}\nonumber\\
&\hspace{0.2cm}-\mathbb{E}\left(\frac{1}{2}R_{1,\ell}(Y)R_{1,j} (Y)\varepsilon_{b_n}^2(Y)+\frac{1}{2}R_{1,\ell}(Y)\varepsilon_{b_n}^2(Y)X_{j} \mathbf{1}_{\{ X_{j} \geqslant 0 \}}\right)\bigg\}+\mathrm{o}_{p}(1)\nonumber\\
&=\frac{1}{\sqrt{n}}\sum\limits_{i=1}^{n}\bigg\{\frac{g_{1,\ell}(Y_i) g_{1,j}(Y_i)}{2f_{b_n}^2(Y_i)}+\frac{g_{1,\ell}(Y_i)f(Y_i)}{2f_{b_n}^2(Y_i)}X_{ij}\mathbf{1}_{\left\{X_{ij}\geqslant 0\right\}}\notag\\
&\hspace{4cm}-\mathbb{E}\left[\frac{g_{1,\ell}(Y) g_{1,j}(Y)}{2f_{b_n}^2(Y)}+\frac{g_{1,\ell}(Y)f(Y)}{2f_{b_n}^2(Y)}X_{j}\mathbf{1}_{\left\{X_{j}\geqslant 0\right\}}\right]\bigg\}+\mathrm{o}_{p}(1).
\end{align}
Then, the proof of $(i)$ is complete by adding to  \eqref{res1i}  the same equation obtained by exchanging $\ell$ and $j$. Since $\widehat{g}^{+}_{1,j,n}$ is obtained from $\widehat{g}^{-}_{1,j,n}$ by replacing $\beta_n$ by $1/{\beta_n}$, we prove $(ii)$ from a similar reasoning. It leads to  
\begin{align*}
&\frac{1}{\sqrt{n}}\sum_{i=1}^n \left\{\frac{g_{1,\ell}(Y_i)}{f_{b_n}^2(Y_i)}\widehat{g}^{+}_{1,j,n}(Y_{i}) - \mathbb{E}\left[ \frac{g_{1,\ell}(Y)}{f_{b_n}^2(Y)}\widehat{g}^{+}_{1,j,n}(Y) \right] \right\}\\
&=\frac{1}{\sqrt{n}}\sum\limits_{i=1}^n\bigg\{\frac{1}{2}R_{1,\ell}(Y_i)R_{1,j} (Y_i)\varepsilon_{b_n}^2(Y_i)+\frac{1}{2}R_{1,\ell}(Y_i)\varepsilon_{b_n}^2(Y_i)X_{ij} \mathbf{1}_{\{ X_{ij} \geqslant 0 \}}\nonumber\\
&\hspace{0.2cm}-\mathbb{E}\left(\frac{1}{2}R_{1,\ell}(Y)R_{1,j} (Y)\varepsilon_{b_n}^2(Y)+\frac{1}{2}R_{1,\ell}(Y)\varepsilon_{b_n}^2(Y)X_{j} \mathbf{1}_{\{ X_{j} \geqslant 0 \}}\right)\bigg\}+\delta_n^\prime+\mathrm{o}_{p}(1),
\end{align*}
where
\begin{align*}
\delta_n^\prime&=\frac{1-\beta_n}{\beta_n\sqrt{n}}\sum\limits_{i=1}^n\bigg\{\frac{1}{2}R_{1,\ell}(Y_i)R_{1,j} (Y_i)\varepsilon_{b_n}^2(Y_i )+\frac{1}{2}R_{1,\ell}(Y_i)\varepsilon_{b_n}^2(Y_i)X_{ij} \mathbf{1}_{\{ X_{ij} \geqslant 0 \}}\nonumber\\
&\hspace{0.2cm}-\mathbb{E}\left(\frac{1}{2}R_{1,\ell}(Y)R_{1,j} (Y)\varepsilon_{b_n}^2(Y)+\frac{1}{2}R_{1,\ell}(Y)\varepsilon_{b_n}^2(Y)X_{j} \mathbf{1}_{\{ X_{j} \geqslant 0 \}}\right)\bigg\}
\end{align*}
and
\begin{align*}
\vert\delta_n^\prime\vert&\leqslant\frac{(1-\beta_n)\sqrt{n}}{2\beta_n}\bigg\{\frac{1}{n}\sum\limits_{i=1}^n\left\vert R_{1,\ell}(Y_i)R_{1,j} (Y_i)\right\vert+\frac{1}{n}\sum\limits_{i=1}^n\left\vert R_{1,\ell}(Y_i)X_{ij} \right\vert+\mathbb{E}\left(\left\vert R_{1,\ell}(Y)R_{1,j} (Y)\right\vert+\left\vert R_{1,\ell}(Y)X_{j} \right\vert\right)\bigg\}.
\end{align*}
Since $(1-\beta_n)\sqrt{n}\sim n^{-7/2}$ and $\beta_n\rightarrow 1$ as $n\rightarrow \infty$, we also deduce  that $\delta_n^\prime=\mathrm{o}_{p}(1)$, what implies the result of $(ii)$ as above.
\hfill $\Box$

\bigskip

\noindent Now,  considering, for any $\ell\in\{1,\ldots,d\}$,
\begin{equation}\label{g1ln}
\widehat{g}_{1,\ell,n}(y)=\frac{1}{nH_{n}(y)}\sum\limits_{i=1}^{n}X_{i\ell}\mathbf{1}_{\left\{X_{i\ell}\geqslant 0\right\}}K\left(\frac{Y_{i}-y}{H_{n}(y)}\right)
\end{equation}
and
\begin{equation}\label{g2ln}
\widehat{g}_{2,\ell,n}(y)=\frac{1}{nH_{n}(y)}\sum\limits_{i=1}^{n}-X_{i\ell}\mathbf{1}_{\left\{X_{i\ell}< 0\right\}}K\left(\frac{Y_{i}-y}{H_{n}(y)}\right),
\end{equation}
we have:

\bigskip

\begin{Lem}\label{l8to11}
Under the assumptions  \ref{as0}, \ref{as3}, \ref{as9} to  \ref{as5} and \ref{as6} , we have for  $(\ell,j)\in\{1,\ldots,d\}^2$ and $(m,r)\in\{1,2\}^2$:
\begin{align*}
&\frac{1}{\sqrt{n}}\sum\limits_{i=1}^{n}\left\{\frac{g_{m,\ell}(Y_i)}{f_{b_n}^2(Y_i)}\widehat{g}_{r,j,n}(Y_{i})+\frac{g_{m,j}(Y_i)}{f_{b_n}^2(Y_i)}\widehat{g}_{r,\ell,n}(Y_{i})-\mathbb{E}\left[\frac{g_{m,\ell}(Y)}{f_{b_n}^2(Y)}\widehat{g}_{r,j,n}(Y)+\frac{g_{m,j}(Y)}{f_{b_n}^2(Y)}\widehat{g}_{r,\ell,n}(Y)\right]\right\}\notag\\
&=\frac{1}{\sqrt{n}}\sum\limits_{i=1}^{n}\bigg\{\frac{g_{m,\ell}(Y_i)g_{r,j}(Y_i)}{f_{b_n}^2(Y_i)}+\frac{g_{m,\ell}(Y_i)f(Y_i)}{2f_{b_n}^2(Y_i)}X_{ij}\mathbf{1}_{\left\{X_{ij}\geqslant 0\right\}}+\frac{g_{m,j}(Y_i)f(Y_i)}{2f_{b_n}^2(Y_i)}X_{i\ell}\mathbf{1}_{\left\{X_{i\ell}\geqslant 0\right\}}\notag\\
&\hspace{0.5cm}-\mathbb{E}\left[\frac{g_{m,\ell}(Y)g_{r,j}(Y)}{f_{b_n}^2(Y)}+\frac{g_{m,\ell}(Y)f(Y)}{2f_{b_n}^2(Y)}X_{j}\mathbf{1}_{\left\{X_{j}\geqslant 0\right\}}+\frac{g_{m,j}(Y)f(Y)}{2f_{b_n}^2(Y)}X_{i\ell}\mathbf{1}_{\left\{X_{\ell}\geqslant 0\right\}} \right]\bigg{\}}+\mathrm{o}_{p}(1).
\end{align*}
\end{Lem}
\noindent\textit{Proof.} We only give the proof for $m=r=1$  since those  related to the other values  are obtained from similar reasoning.
By using Assumption \ref{as4}-$(vi)$  we have, as in the proof of Theorem 2 in \cite{bengono24} (see p. 14), $\widehat{g}_{1,\ell,n}^{-}(y)\leqslant\widehat{g}_{1,\ell,n}(y)\leqslant \widehat{g}_{1,\ell,n}^{+}(y) $, where $\widehat{g}_{1,\ell,n}^{-}$ and $\widehat{g}_{1,\ell,n}^{+}$ are defined in \eqref{g1m} and \eqref{g1p}. Hence
\begin{align*}
&\frac{1}{\sqrt{n}}\sum\limits_{i=1}^{n}\left\{\frac{g_{1,\ell}(Y_i)}{f_{b_n}^2(Y_i)}\widehat{g}^{-}_{1,j,n}(Y_{i})-\mathbb{E}\left[\frac{g_{1,\ell}(Y)}{f_{b_n}^2(Y)}\widehat{g}_{1,n,j}^{+}(Y)\right]\right\}\notag\\
&\leqslant\frac{1}{\sqrt{n}}\sum\limits_{i=1}^{n}\left\{\frac{g_{1,\ell}(Y_i)}{f_{b_n}^2(Y_i)}\widehat{g}_{1,j,n}(Y_{i})-\mathbb{E}\left[\frac{g_{1,\ell}(Y)}{f_{b_n}^2(Y)}\widehat{g}_{1,j,n}(Y)\right]\right\}
\leqslant\frac{1}{\sqrt{n}}\sum\limits_{i=1}^{n}\left\{\frac{g_{1,\ell}(Y_i)}{f_{b_n}^2(Y_i)}\widehat{g}^{+}_{1,j,n}(Y_{i})-\mathbb{E}\left[\frac{g_{1,\ell}(Y)}{f_{b_n}^2(Y)}\widehat{g}_{1,j,n}^{-}(Y)\right]\right\}.
\end{align*}
However, by using Lemma \ref{l3}, we obtain
\begin{align*}
&\frac{1}{\sqrt{n}}\sum\limits_{i=1}^{n}\left\{\frac{g_{1,\ell}(Y_i)}{f_{b_n}^2(Y_i)}\widehat{g}^{-}_{1,j,n}(Y_{i})-\mathbb{E}\left[\frac{g_{1,\ell}(Y)}{f_{b_n}^2(Y)}\widehat{g}_{1,n,j}^{+}(Y)\right]\right\}\\
&=\frac{1}{\sqrt{n}}\sum\limits_{i=1}^{n}\left\{\frac{g_{1,\ell}(Y_i)}{f_{b_n}^2(Y_i)}\widehat{g}^{-}_{1,j,n}(Y_{i})-\mathbb{E}\left[\frac{g_{1,\ell}(Y)}{f_{b_n}^2(Y)}\widehat{g}_{1,n,j}^{-}(Y)\right]\right\}+ \sqrt{n}\bigg(\mathbb{E}\bigg[\frac{g_{1,\ell}(Y)}{f_{b_n}^2(Y)}\Big(\widehat{g}^{-}_{1,j,n}(Y)-\widehat{g}^{+}_{1,j,n}(Y)\Big)\bigg]\bigg)\\
&=\frac{1}{\sqrt{n}}\sum\limits_{i=1}^{n}\left\{\frac{g_{1,\ell}(Y_i)}{f_{b_n}^2(Y_i)}\widehat{g}^{-}_{1,j,n}(Y_{i})-\mathbb{E}\left[\frac{g_{1,\ell}(Y)}{f_{b_n}^2(Y)}\widehat{g}_{1,n,j}^{-}(Y)\right]\right\}+\mathrm{o}(1)\\
\end{align*}
and, similarly,
\begin{align*}
&\frac{1}{\sqrt{n}}\sum\limits_{i=1}^{n}\left\{\frac{g_{1,\ell}(Y_i)}{f_{b_n}^2(Y_i)}\widehat{g}^{+}_{1,j,n}(Y_{i})-\mathbb{E}\left[\frac{g_{1,\ell}(Y)}{f_{b_n}^2(Y)}\widehat{g}_{1,j,n}^{-}(Y)\right]\right\}\\
&=\frac{1}{\sqrt{n}}\sum\limits_{i=1}^{n}\left\{\frac{g_{1,\ell}(Y_i)}{f_{b_n}^2(Y_i)}\widehat{g}^{+}_{1,j,n}(Y_{i})-\mathbb{E}\left[\frac{g_{1,\ell}(Y)}{f_{b_n}^2(Y)}\widehat{g}_{1,j,n}^{+}(Y)\right]\right\}+\mathrm{o}(1).
\end{align*}
Hence
\begin{align*}
&\frac{1}{\sqrt{n}}\sum\limits_{i=1}^{n}\left\{\frac{g_{1,\ell}(Y_i)}{f_{b_n}^2(Y_i)}\widehat{g}^{-}_{1,j,n}(Y_{i})-\mathbb{E}\left[\frac{g_{1,\ell}(Y)}{f_{b_n}^2(Y)}\widehat{g}_{1,j,n}^{-}(Y)\right]\right\}+\mathrm{o}(1)\\
&\hspace{2cm}\leqslant\frac{1}{\sqrt{n}}\sum\limits_{i=1}^{n}\left\{\frac{g_{1,\ell}(Y_i)}{f_{b_n}^2(Y_i)}\widehat{g}_{1,j,n}(Y_{i})-\mathbb{E}\left[\frac{g_{1,\ell}(Y)}{f_{b_n}^2(Y)}\widehat{g}_{1,j,n}(Y)\right]\right\}\\
&\hspace{5cm}\leqslant\frac{1}{\sqrt{n}}\sum\limits_{i=1}^{n}\left\{\frac{g_{1,\ell}(Y_i)}{f_{b_n}^2(Y_i)}\widehat{g}^{+}_{1,j,n}(Y_{i})-\mathbb{E}\left[\frac{g_{1,\ell}(Y)}{f_{b_n}^2(Y)}\widehat{g}_{1,j,n}^{+}(Y)\right]\right\}+\mathrm{o}(1).
\end{align*}
By adding to this last inequality   the same one obtained by exchanging $\ell$ and $j$, we obtain $ \mathscr{F}_{\ell,j,n}^{-}+\mathrm{o}(1)\leqslant \mathscr{F}_{\ell,j,n}\leqslant  \mathscr{F}_{\ell,j,n}^{+}+\mathrm{o}(1)$, where
\begin{align*}
 & \mathscr{F}_{\ell,j,n}=\frac{1}{\sqrt{n}}\sum\limits_{i=1}^{n}\left\{\frac{g_{1,\ell}(Y_i)}{f_{b_n}^2(Y_i)}\widehat{g}_{1,j,n}(Y_{i})+\frac{g_{1,j}(Y_i)}{f_{b_n}^2(Y_i)}\widehat{g}_{1,\ell,n}(Y_{i})-\mathbb{E}\left[\frac{g_{1,\ell}(Y)}{f_{b_n}^2(Y)}\widehat{g}_{1,j,n}(Y)+\frac{g_{1,j}(Y)}{f_{b_n}^2(Y)}\widehat{g}_{1,\ell,n}(Y)\right]\right\},\\
& \mathscr{F}_{\ell,j,n}^{-}=\frac{1}{\sqrt{n}}\sum\limits_{i=1}^{n}\left\{\frac{g_{1,\ell}(Y_i)}{f_{b_n}^2(Y_i)}\widehat{g}_{1,j,n}^{-}(Y_{i})+\frac{g_{1,j}(Y_i)}{f_{b_n}^2(Y_i)}\widehat{g}_{1,\ell,n}^{-}(Y_{i})-\mathbb{E}\left[\frac{g_{1,\ell}(Y)}{f_{b_n}^2(Y)}\widehat{g}_{1,j,n}^{-}(Y)+\frac{g_{1,j}(Y)}{f_{b_n}^2(Y)}\widehat{g}_{1,\ell,n}^{-}(Y)\right]\right\},\\
& \mathscr{F}_{\ell,j,n}^{+} =\frac{1}{\sqrt{n}}\sum\limits_{i=1}^{n}\left\{\frac{g_{1,\ell}(Y_i)}{f_{b_n}^2(Y_i)}\widehat{g}_{1,j,n}^{+}(Y_{i})+\frac{g_{1,j}(Y_i)}{f_{b_n}^2(Y_i)}\widehat{g}_{1,\ell,n}^{+}(Y_{i})-\mathbb{E}\left[\frac{g_{1,\ell}(Y)}{f_{b_n}^2(Y)}\widehat{g}_{1,j,n}^{+}(Y)+\frac{g_{1,j}(Y)}{f_{b_n}^2(Y)}\widehat{g}_{1,\ell,n}^{+}(Y)\right]\right\}.
\end{align*}
Then, by using Lemma \ref{l67}, we deduce from the preceding inequality that
\begin{align*}
& \mathscr{F}_{\ell,j,n}-\frac{1}{\sqrt{n}}\sum\limits_{i=1}^{n}\bigg\{\frac{g_{1,\ell}(Y_i)g_{1,j}(Y_i)}{f_{b_n}^2(Y_i)}+\frac{g_{1,\ell}(Y_i)f(Y_i)}{2f_{b_n}^2(Y_i)}X_{ij}\mathbf{1}_{\left\{X_{ij}\geqslant 0\right\}}+\frac{g_{1,j}(Y_i)f(Y_i)}{2f_{b_n}^2(Y_i)}X_{i\ell}\mathbf{1}_{\left\{X_{i\ell}\geqslant 0\right\}}\\
&+\mathbb{E}\left[\frac{g_{1,\ell}(Y)g_{1,j}(Y)}{f_{b_n}^2(Y)}+\frac{g_{1,\ell}(Y)f(Y)}{2f_{b_n}^2(Y_i)}X_{j}\mathbf{1}_{\left\{X_{j}\geqslant 0\right\}}+\frac{g_{1,j}(Y)f(Y)}{2f_{b_n}^2(Y)}X_{\ell}\mathbf{1}_{\left\{X_{\ell}\geqslant 0\right\}} \right]\bigg\}=\mathrm{o}_{p}(1),
\end{align*}
and the proof  is complete.
\hfill $\Box$

\bigskip

\noindent We define the functions
\begin{equation}\label{Ilj2}
I_{\ell j}^{(2)}(y) =\frac{g_\ell(y)\widehat{g}_{j , n}(y)+g_{j}(y)\widehat{g}_{\ell , n}(y)}{f_{b_n}^2(y)}=\frac{R_{b_n,\ell}(y)\widehat{g}_{j , n}(y)}{f_{b_n}(y)} + \frac{R_{b_n,j}(y)\widehat{g}_{\ell , n}(y)}{f_{b_n}(y)}
\end{equation}
and
\begin{equation}\label{Ilj3}
I_{\ell j}^{(3)}(y) = 2R_{b_n,\ell}(y)R_{b_n,j}(y)\frac{\widehat{f}_{b_n}(y)}{{f}_{b_n}(y)},
\end{equation}
where $R_{b_n,\ell}$ is defined in \eqref{rbnl}. Then, we have:

\bigskip

\begin{Lem}\label{l12to13}
Under  the assumptions  \ref{as0} and  \ref{as1} to  \ref{as6}, we have for  $(\ell,j)\in\{1,\ldots,d\}^2$:
\begin{align*}
(i)\,\,&\frac{1}{\sqrt{n}}\sum\limits_{i=1}^{n}\left\{I_{\ell j}^{(2)}(Y_{i})-\mathbb{E}\left[I_{\ell j}^{(2)}(Y)\right]\right\}\nonumber\\
&=\frac{1}{\sqrt{n}}\sum_{i=1}^n\bigg\{R_{b_n,\ell}\left(Y_i\right)R_{b_n,j}\left(Y_i\right)+ \frac{1}{2}X_{ij}R_{b_n,\ell}(Y_i)\frac{f(Y_i)}{f_{b_n}(Y_i)}+\frac{1}{2}X_{i\ell}R_{b_n,j}(Y_i)\frac{f(Y_i)}{f_{b_n}(Y_i)}\\
&\hspace{0.5cm} -\mathbb{E} \left[ R_{b_n,\ell}\left(Y\right)R_{b_n,j}\left(Y\right)+ \frac{1}{2}X_{j}R_{b_n,\ell}(Y)\frac{f(Y)}{f_{b_n}(Y)}+\frac{1}{2}X_{\ell}R_{b_n,j}(Y)\frac{f(Y)}{f_{b_n}(Y)} \right]\bigg\}+ \mathrm{o}_{p}(1);
\end{align*}
\begin{align*}
(ii)\,\,&\frac{1}{\sqrt{n}}\sum\limits_{i=1}^{n}\left\{I_{\ell j}^{(3)}(Y_{i})-\mathbb{E}\left[I_{\ell j}^{(3)}(Y)\right]\right\}\nonumber\\
&=\frac{2}{\sqrt{n}}\sum\limits_{i=1}^n\left\{R_{b_n,\ell}\left(Y_i\right)R_{b_n,j}\left(Y_i\right)\frac{f(Y_{i})}{f_{b_{n}}(Y_{i})}-\mathbb{E}\left[R_{b_n,\ell}\left(Y\right)R_{b_n,j}\left(Y\right)\frac{f(Y)}{f_{b_{n}}(Y)}\right]\right\}+\mathrm{o}_{p}(1).
\end{align*}
\end{Lem}
\noindent\textit{Proof.} 

\noindent\textit{(i)}. We have:
\begin{align*}
&\frac{1}{\sqrt{n}}\sum\limits_{i=1}^{n}\left\{I_{\ell j}^{(2)}(Y_{i})-\mathbb{E}\left[I_{\ell j}^{(2)}(Y)\right]\right\}\nonumber\\
&=\frac{1}{\sqrt{n}}\sum_{i=1}^n\bigg\{\frac{g_{\ell}(Y_i)\widehat{g}_{j , n}(Y_i)}{f_{b_n}^2(Y_i)} + \frac{g_{j}(Y_i)\widehat{g}_{\ell , n}(Y_i)}{f_{b_n}^2(Y_i)} -\mathbb{E} \left[ \frac{g_{\ell}(Y)\widehat{g}_{j , n}(Y)}{f_{b_n}^2(Y)} + \frac{g_{j}(Y)\widehat{g}_{\ell , n}(Y)}{f_{b_n}^2(Y)} \right]\bigg\}.
\end{align*}
Since $g_{\ell}=g_{1,\ell}-g_{2,\ell}$ and $\widehat{g}_{\ell, n}=\widehat{g}_{1,\ell, n}-\widehat{g}_{2,\ell , n}$, where $\widehat{g}_{1,\ell, n}$ and $\widehat{g}_{2,\ell, n}$ are defined in \eqref{g1ln} and  \eqref{g2ln},  it follows
\begin{align*}
&\frac{1}{\sqrt{n}}\sum\limits_{i=1}^{n}\left\{I_{\ell j}^{(2)}(Y_{i})-\mathbb{E}\left[I_{\ell j}^{(2)}(Y)\right]\right\}\nonumber\\
&=\frac{1}{\sqrt{n}}\sum\limits_{i=1}^{n}\left\{\frac{g_{1,\ell}(Y_i)}{f_{b_n}^2(Y_i)}\widehat{g}_{1,j,n}(Y_{i})+\frac{g_{1,j}(Y_i)}{f_{b_n}^2(Y_i)}\widehat{g}_{1,\ell,n}(Y_{i})-\mathbb{E}\left[\frac{g_{1,\ell}(Y)}{f_{b_n}^2(Y)}\widehat{g}_{1,j,n}(Y)+\frac{g_{1,j}(Y)}{f_{b_n}^2(Y)}\widehat{g}_{1,\ell,n}(Y)\right]\right\}\\
&-\frac{1}{\sqrt{n}}\sum\limits_{i=1}^{n}\left\{\frac{g_{1,\ell}(Y_i)}{f_{b_n}^2(Y_i)}\widehat{g}_{2,j,n}(Y_{i})+\frac{g_{1,j}(Y_i)}{f_{b_n}^2(Y_i)}\widehat{g}_{2,\ell,n}(Y_{i})-\mathbb{E}\left[\frac{g_{1,\ell}(Y)}{f_{b_n}^2(Y)}\widehat{g}_{2,j,n}(Y)+\frac{g_{1,j}(Y)}{f_{b_n}^2(Y)}\widehat{g}_{2,\ell,n}(Y)\right]\right\}\\
&-\frac{1}{\sqrt{n}}\sum\limits_{i=1}^{n}\left\{\frac{g_{2,\ell}(Y_i)}{f_{b_n}^2(Y_i)}\widehat{g}_{1,j,n}(Y_{i})+\frac{g_{2,j}(Y_i)}{f_{b_n}^2(Y_i)}\widehat{g}_{1,\ell,n}(Y_{i})-\mathbb{E}\left[\frac{g_{2,\ell}(Y)}{f_{b_n}^2(Y)}\widehat{g}_{1,j,n}(Y)+\frac{g_{2,j}(Y)}{f_{b_n}^2(Y)}\widehat{g}_{1,\ell,n}(Y)\right]\right\}\\
&+\frac{1}{\sqrt{n}}\sum\limits_{i=1}^{n}\left\{\frac{g_{2,\ell}(Y_i)}{f_{b_n}^2(Y_i)}\widehat{g}_{2,j,n}(Y_{i})+\frac{g_{2,j}(Y_i)}{f_{b_n}^2(Y_i)}\widehat{g}_{2,\ell,n}(Y_{i})-\mathbb{E}\left[\frac{g_{2,\ell}(Y)}{f_{b_n}^2(Y)}\widehat{g}_{2,j,n}(Y)+\frac{g_{2,j}(Y)}{f_{b_n}^2(Y)}\widehat{g}_{2,\ell,n}(Y)\right]\right\}.
\end{align*}
Then, Lemma \ref{l8to11} yields the result.

\noindent\textit{(ii)}. Putting
\[
\mathcal{Z}_{\ell,j,n}=\frac{2}{\sqrt{n}}\sum\limits_{i=1}^n\left\{R_{b_n,\ell}\left(Y_i\right)R_{b_n,j}\left(Y_i\right)\frac{f(Y_{i})}{f_{b_{n}}(Y_{i})}-\mathbb{E}\left[R_{b_n,\ell}\left(Y\right)R_{b_n,j}\left(Y\right)\frac{f(Y)}{f_{b_{n}}(Y)}\right]\right\},
\]
we have:
\begin{align*}
&\mathbb{E}\bigg[\bigg(\frac{1}{\sqrt{n}}\sum\limits_{i=1}^{n}\left\{I_{\ell j}^{(3)}(Y_{i})-\mathbb{E}\left[I_{\ell j}^{(3)}(Y)\right]\right\}-\mathcal{Z}_{\ell,j,n}\bigg)^2\bigg]\\
&=\mathbb{E}\bigg[\bigg(\frac{2}{\sqrt{n}}\sum\limits_{i=1}^n\bigg\{R_{b_n,\ell}\left(Y_i\right)R_{b_n,j}\left(Y_i\right)\frac{\widehat{f}_{n}(Y_i)-f(Y_{i})}{f_{b_{n}}(Y_{i})}-\mathbb{E}\left[R_{b_n,\ell}\left(Y\right)R_{b_n,j}\left(Y\right)\frac{\widehat{f}_{n}(Y)-f(Y)}{f_{b_{n}}(Y)}\right]\bigg\}\bigg)^{2}\bigg]\\
&\leqslant 4\mathbb{E}\left[\left(R_{b_n,\ell}\left(Y\right)R_{b_n,j}\left(Y\right)\frac{\widehat{f}_{n}(Y)-f(Y)}{f_{b_{n}}(Y)}\right)^2\right].
\end{align*}
Using Theorem 1 of \cite{bengono24} together with the inequalities $f_{b_{n}}\geqslant b_n$ and $R_{b_n,\ell}^2=R_{\ell}^2\varepsilon_{b_n}^2\leqslant R_{\ell}^2$, we obtain almost surely
\[
\mathbb{E}\bigg[\bigg(\frac{1}{\sqrt{n}}\sum\limits_{i=1}^{n}\left\{I_{\ell j}^{(3)}(Y_{i})-\mathbb{E}\left[I_{\ell j}^{(3)}(Y)\right]\right\}-\mathcal{Z}_{\ell,j,n}\bigg)^2\bigg]\leqslant 4b_n^{-2}\tau_n^2\mathbb{E}\left[R_{\ell}^2\left(Y\right)R_{j}^2\left(Y\right)\right],
\]
where $\tau_n=\frac{k_{n}^4}{n^4}+\frac{\sqrt{n\log(n)}}{k_{n}}$. From \eqref{equivrt}, we have  $\tau_n\sim n^{1/2-c_1}\log^{1/2}(n)$ and, therefore,    $b_n^{-2}\tau_n^2\sim  n^{2\left(1/2-c_1+c_2\right)}\log(n)$. Since $c_1-c_2>3/4>1/2$, it follows that  $b_n^{-2}\tau_n^2\rightarrow 0$ as $n\rightarrow\infty$ and, consequently, that
\[
\mathbb{E}\bigg[\bigg(\frac{1}{\sqrt{n}}\sum\limits_{i=1}^{n}\left\{I_{\ell j}^{(3)}(Y_{i})-\mathbb{E}\left[I_{\ell j}^{(3)}(Y)\right]\right\}-\mathcal{Z}_{\ell,j,n}\bigg)^2\bigg]=\mathrm{o}(1).
\]
Then, using Bienaym\'e-Chebyshev inequality we deduce that
\[
\frac{1}{\sqrt{n}}\sum\limits_{i=1}^{n}\left\{I_{\ell j}^{(3)}(Y_{i})-\mathbb{E}\left[I_{\ell j}^{(3)}(Y)\right]\right\}-\mathcal{Z}_{\ell,j,n}=\mathrm{o}_p(1),
\]
and the proof of $(ii)$ is complete.
\hfill $\Box$ 

\bigskip

\begin{Lem}\label{l14}
Under the assumptions  \ref{as3} and \ref{as9} to \ref{as6}, we have for  $(\ell,j)\in\{1,\ldots,d\}^2$:
\begin{equation*}
\sqrt{n}\mathbb{E}\left[\frac{R_{b_n,\ell}\left(Y\right)\widehat{g}_{j,n}(Y)}{f_{b_{n}}(Y)}\right]=\sqrt{n}\mathbb{E}\left[R_{\ell}\left(Y\right)R_{j}\left(Y\right)\right]+\mathrm{o}(1).
\end{equation*}
\end{Lem}
\noindent\textit{Proof.} 
Clearly
\begin{align}\label{decl14}
&\frac{R_{b_n,\ell}\left(Y\right)\widehat{g}_{j,n}(Y)}{f_{b_{n}}(Y)}\nonumber\\
&=\frac{g_\ell\left(Y\right)\widehat{g}_{j,n}(Y)}{f_{b_{n}}^2(Y)}=\frac{\left(g_{1,\ell}\left(Y\right)-g_{2,\ell}\left(Y\right)\right)\left(\widehat{g}_{1,j,n}(Y)-\widehat{g}_{2,j,n}(Y)\right)}{f_{b_{n}}^2(Y)}\nonumber\\
&=\frac{g_{1,\ell}(Y)\widehat{g}_{1,j,n}(Y)}{f_{b_n}^2(Y)}-\frac{g_{1,\ell}(Y)\widehat{g}_{2,j,n}(Y)}{f_{b_n}^2(Y)}-\frac{g_{2,\ell}(Y)\widehat{g}_{1,j,n}(Y)}{f_{b_n}^2(Y)}+\frac{g_{2,\ell}(Y)\widehat{g}_{2,j,n}(Y)}{f_{b_n}^2(Y)}.
\end{align}
Using  $\widehat{g}_{1,\ell,n}^{-}(y)\leqslant\widehat{g}_{1,\ell,n}(y)\leqslant \widehat{g}_{1,\ell,n}^{+}(y) $ (see \cite{bengono24},  p. 14), we obtain
\[
\sqrt{n}\mathbb{E}\left[\frac{g_{1,\ell}(Y)\widehat{g}_{1,\ell,n}^{-}(Y)}{f_{b_n}^2(Y)}\right]\leqslant \sqrt{n}\mathbb{E}\left[\frac{g_{1,\ell}(Y)\widehat{g}_{1,\ell,n}(Y)}{f_{b_n}^2(Y)}\right]\leqslant \sqrt{n}\mathbb{E}\left[\frac{g_{1,\ell}(Y)\widehat{g}_{1,\ell,n}^{+}(Y)}{f_{b_n}^2(Y)}\right];
\]
then, \eqref{res1} and \eqref{res2} lead to
\[
\mathrm{o}(1)\leqslant \sqrt{n}\mathbb{E}\left[\frac{g_{1,\ell}(Y)\widehat{g}_{1,\ell,n}(Y)}{f_{b_n}^2(Y)}\right]-\sqrt{n}\mathbb{E}\left[R_{1,\ell}\left(Y\right)R_{1,j}\left(Y\right)\right]\leqslant \mathrm{o}(1),
\]
what allows to conclude that 
\begin{equation*}
\sqrt{n}\mathbb{E}\left[\frac{g_{1,\ell}(Y)\widehat{g}_{1,\ell,n}(Y)}{f_{b_n}^2(Y)}\right]=\sqrt{n}\mathbb{E}\left[R_{1,\ell}\left(Y\right)R_{1,j}\left(Y\right)\right]+\mathrm{o}(1).
\end{equation*}
From a similar reasoning, we also obtain
\begin{equation*}
\sqrt{n}\mathbb{E}\left[\frac{g_{1,\ell}(Y)\widehat{g}_{2,j,n}(Y)}{f_{b_n}^2(Y)}\right]=\sqrt{n}\mathbb{E}\left[R_{1,\ell}\left(Y\right)R_{2,j}\left(Y\right)\right]+\mathrm{o}(1),
\end{equation*}
\begin{equation*}
\sqrt{n}\mathbb{E}\left[\frac{g_{2,\ell}(Y)\widehat{g}_{1,j,n}(Y)}{f_{b_n}^2(Y)}\right]=\sqrt{n}\mathbb{E}\left[R_{2,\ell}\left(Y\right)R_{1,j}\left(Y\right)\right]+\mathrm{o}(1),
\end{equation*}
and
\begin{equation*}
\sqrt{n}\mathbb{E}\left[\frac{g_{2,\ell}(Y)\widehat{g}_{2,j,n}(Y)}{f_{b_n}^2(Y)}\right]=\sqrt{n}\mathbb{E}\left[R_{2,\ell}\left(Y\right)R_{2,j}\left(Y\right)\right]+\mathrm{o}(1).
\end{equation*}
Then, \eqref{decl14} and the decomposition
\begin{equation}\label{decR}
R_{\ell}\left(Y\right)R_{j}\left(Y\right)=R_{1,\ell}\left(Y\right)R_{1,j}\left(Y\right)-R_{1,\ell}\left(Y\right)R_{2,j}\left(Y\right)-R_{2,\ell}\left(Y\right)R_{1,j}\left(Y\right)+R_{2,\ell}\left(Y\right)R_{2,j}\left(Y\right)
\end{equation}
yield  the result.
\hfill $\Box$

\bigskip

\noindent Let us put
 \begin{equation}\label{f1n}
\widehat{f}_{1,n}(y)=\frac{1}{nD_{n}^{+}(y)}\sum\limits_{i=1}^{n}K\left(\frac{Y_{i}-y}{D_{n}^{-}(y)}\right)
\end{equation}
and
\begin{equation}\label{f2n}
\widehat{f}_{2,n}(y)=\frac{1}{nD_{n}^{-}(y)}\sum\limits_{i=1}^{n}K\left(\frac{Y_{i}-y}{D_{n}^{+}(y)}\right),
\end{equation}
where $D_{n}^{-}$ and $D_{n}^{+}$ are defined in \eqref{dnmp}. Then, we have:

\bigskip

\begin{Lem}\label{l15}
Under the assumptions \ref{as2}, \ref{as8}  and \ref{as9}, we have for  $(\ell,j)\in\{1,\ldots,d\}^2$,  $(m,r,s)\in\{1,2\}^3$ and any $C>0$:
\[
\sqrt{n}\mathbb{E}\left[\frac{R_{m,\ell}\left(Y\right)R_{s,j}\left(Y\right)\widehat{f}_{r,n}(Y)}{f(Y)}\mathbf{1}_{\left\{f(Y)\geqslant b_{n}+C\tau_{n}\right\}}\right]=\sqrt{n}\mathbb{E}\left[R_{m,\ell}\left(Y\right)R_{s,j}\left(Y\right)\right]+\mathrm{o}(1).
\]
\end{Lem}
\noindent\textit{Proof.} We have
\begin{align*}
&\sqrt{n}\mathbb{E}\left[\frac{R_{m,\ell}\left(Y\right)R_{s,j}\left(Y\right)\widehat{f}_{1,n}(Y)}{f(Y)}\mathbf{1}_{\left\{f(Y)\geqslant b_{n}+C\tau_{n}\right\}}\right]\\
&=\frac{\sqrt{n\beta_n}}{k_n}\sum_{i=1}^n\mathbb{E}\left[R_{m,\ell}\left(Y\right)R_{s,j}\left(Y\right)K\bigg(\frac{nf(Y)(Y_i-Y)}{k_n\sqrt{\beta_n}}\bigg)\mathbf{1}_{\left\{f(Y)\geqslant b_{n}+C\tau_{n}\right\}}\right]\\
&=\frac{n^{3/2}\sqrt{\beta_n}}{k_n} \mathbb{E}\left[R_{m,\ell}\left(Y\right)R_{s,j}\left(Y\right)K\bigg(\frac{nf(Y)(Y_1-Y)}{k_n\sqrt{\beta_n}}\bigg)\mathbf{1}_{\left\{f(Y)\geqslant b_{n}+C\tau_{n}\right\}}\right]\\
&=\frac{n^{3/2}\sqrt{\beta_n}}{k_n} \int\int_{\left\{f(y)\geqslant b_{n}+C\tau_{n}\right\}} R_{m,\ell}\left(y\right)R_{s,j}\left(y\right)K\bigg(\frac{nf(y)(z-y)}{k_n\sqrt{\beta_n}}\bigg)f(y)\,f(z)\,\,dy\,dz.
\end{align*}
Since
\[
 \int K\bigg(\frac{nf(y)(z-y)}{k_n\sqrt{\beta_n}}\bigg)\,f(z)\,\,dz=\frac{k_n\sqrt{\beta_n}}{nf(y)}\int f\bigg(y+\frac{k_n\sqrt{\beta_n}}{nf(y)}t\bigg)\,K(t)\,\,dt,
\]
it follows
\begin{align}\label{decl15}
&\sqrt{n}\mathbb{E}\left[\frac{R_{m,\ell}\left(Y\right)R_{s,j}\left(Y\right)\widehat{f}_{1,n}(Y)}{f(Y)}\mathbf{1}_{\left\{f(Y)\geqslant b_{n}+C\tau_{n}\right\}}\right]\nonumber\\
&=  \beta_n\sqrt{n}  \int\int_{\left\{f(y)\geqslant b_{n}+C\tau_{n}\right\}} R_{m,\ell}\left(y\right)R_{s,j}\left(y\right) f\bigg(y+\frac{k_n\sqrt{\beta_n}}{nf(y)}t\bigg)\,K(t)\,\,dy\,dt\nonumber\\
&=  \beta_n\sqrt{n}  \int\int_{\left\{f(y)\geqslant b_{n}+C\tau_{n}\right\}} R_{m,\ell}\left(y\right)R_{s,j}\left(y\right) \bigg(f\left(y+\frac{k_n\sqrt{\beta_n}}{nf(y)}t\right)-f(y)\bigg)\,K(t)\,\,dy\,dt\nonumber\\
&\hspace{0.5cm}+\beta_n\sqrt{n} \int_{\left\{f(y)\geqslant b_{n}+C\tau_{n}\right\}} R_{m,\ell}\left(y\right)R_{s,j}\left(y\right)  f(y)\,\,dy.
\end{align}
However,
\begin{align*}
&\beta_n\sqrt{n} \int_{\left\{f(y)\geqslant b_{n}+C\tau_{n}\right\}} R_{m,\ell}\left(y\right)R_{s,j}\left(y\right)  f(y)\,\,dy\\
&=\beta_{n}\sqrt{n}\mathbb{E}\left[R_{m,\ell}\left(Y\right)R_{s,j}\left(Y\right)\mathbf{1}_{\left\{f(Y)\geqslant b_{n}+C\tau_{n}\right\}}\right]\notag\\
&=\beta_{n}\sqrt{n}\mathbb{E}\left[R_{m,\ell}\left(Y\right)R_{s,j}\left(Y\right)\right]-\beta_{n}\sqrt{n}\mathbb{E}\left[R_{m,\ell}\left(Y\right)R_{s,j}\left(Y\right)\mathbf{1}_{\left\{f(Y)< b_{n}+C\tau_{n}\right\}}\right],
\end{align*}
so that
\begin{align*}
&\beta_n\sqrt{n} \int_{\left\{f(y)\geqslant b_{n}+C\tau_{n}\right\}} R_{m,\ell}\left(y\right)R_{s,j}\left(y\right)  f(y)\,\,dy-\sqrt{n}\mathbb{E}\left[R_{m,\ell}\left(Y\right)R_{s,j}\left(Y\right)\right]\nonumber\\
&=-(1-\beta_{n})\sqrt{n}\mathbb{E}\left[R_{m,\ell}\left(Y\right)R_{s,j}\left(Y\right)\right]-\beta_{n}\sqrt{n}\mathbb{E}\left[R_{m,\ell}\left(Y\right)R_{s,j}\left(Y\right)\mathbf{1}_{\left\{f(Y)< b_{n}+C\tau_{n}\right\}}\right].
\end{align*}
Since $(1-\beta_{n})\sqrt{n}\sim n^{-7/2}$, $\beta_n\rightarrow 1$ as  $n\rightarrow \infty$ and $ b_{n}+C\tau_{n}\sim b_n$ because $b_n^{-1}\tau_n\sim n^{1/2-c_1+c_2}\rightarrow 0$ as $n\rightarrow \infty$, it follows from the above equality and Assumption \ref{as9} that 
\begin{equation}\label{cv1l15}
\beta_n\sqrt{n} \int_{\left\{f(y)\geqslant b_{n}+C\tau_{n}\right\}} R_{m,\ell}\left(y\right)R_{s,j}\left(y\right)  f(y)\,\,dy=\sqrt{n}\mathbb{E}\left[R_{m,\ell}\left(Y\right)R_{s,j}\left(Y\right)\right]+\mathrm{o}(1).
\end{equation}
Moreover, an use of Taylor's  expansion up to order 3  leads to 
\begin{align*}
&\beta_n\sqrt{n}  \int\int_{\left\{f(y)\geqslant b_{n}+C\tau_{n}\right\}} R_{m,\ell}\left(y\right)R_{s,j}\left(y\right) \bigg(f\left(y+\frac{k_n\sqrt{\beta_n}}{nf(y)}t\right)-f(y)\bigg)\,K(t)\,\,dy\,dt\nonumber\\
&=\beta_n\sqrt{n}  \sum_{q=1}^2\int_{\left\{f(y)\geqslant b_{n}+C\tau_{n}\right\}} R_{m,\ell}\left(y\right)R_{s,j}\left(y\right) \frac{f^{(q)}(y)}{q!}\frac{k_n^q(\sqrt{\beta_n})^q}{n^qf^q(y)}\bigg(\int t^{q}\,K(t)\,\,dt\bigg)dy\\
&+\beta_n\sqrt{n}   \int_{\left\{f(y)\geqslant b_{n}+C\tau_{n}\right\}} R_{m,\ell}\left(y\right)R_{s,j}\left(y\right)\frac{k_n^3(\sqrt{\beta_n})^3}{6n^3f^3(y)}\bigg(\int t^{3}f^{(3)}\left(y+\theta\frac{k_n\sqrt{\beta_n}}{nf(y)}t\right)\,K(t)\,\,dt\bigg)dy,
\end{align*}
where $0<\theta<1$. Since $K$ is of order 3, it follows
\begin{align*}
&\beta_n\sqrt{n}  \int\int_{\left\{f(y)\geqslant b_{n}+C\tau_{n}\right\}} R_{m,\ell}\left(y\right)R_{s,j}\left(y\right) \bigg(f\left(y+\frac{k_n\sqrt{\beta_n}}{nf(y)}t\right)-f(y)\bigg)\,K(t)\,\,dy\,dt\nonumber\\
&=\beta_n\sqrt{n}   \int_{\left\{f(y)\geqslant b_{n}+C\tau_{n}\right\}} R_{m,\ell}\left(y\right)R_{s,j}\left(y\right)\frac{k_n^3(\sqrt{\beta_n})^3}{6n^3f^3(y)}\bigg(\int t^{3}\bigg(f^{(3)}\left(y+\theta\frac{k_n\sqrt{\beta_n}}{nf(y)}t\right)-f^{(3)}(y)\bigg)\,K(t)\,\,dt\bigg)dy.
\end{align*}
Thus
\begin{align*}
&\bigg\vert\beta_n\sqrt{n}  \int\int_{\left\{f(y)\geqslant b_{n}+C\tau_{n}\right\}} R_{m,\ell}\left(y\right)R_{s,j}\left(y\right) \bigg(f\left(y+\frac{k_n\sqrt{\beta_n}}{nf(y)}t\right)-f(y)\bigg)\,K(t)\,\,dy\,dt\bigg\vert\\
&\leqslant\beta_n\sqrt{n}   \int_{\left\{f(y)\geqslant b_{n}+C\tau_{n}\right\}} \left\vert R_{m,\ell}\left(y\right)R_{s,j}\left(y\right)\right\vert\frac{k_n^3(\sqrt{\beta_n})^3}{6n^3f^3(y)}\bigg(\int \vert t\vert^3\,\bigg\vert f^{(3)}\left(y+\theta\frac{k_n\sqrt{\beta_n}}{nf(y)}t\right)-f^{(3)}(y)\bigg\vert\,\vert K(t)\vert   \,\,dt\bigg)dy\\
&\leqslant\frac{\beta_n^3n^{-7/2} k_n^4}{6}   \int_{\left\{f(y)\geqslant b_{n}+C\tau_{n}\right\}}  \frac{ \left\vert R_{m,\ell}\left(y\right)R_{s,j}\left(y\right)\right\vert}{f^4(y)}\,\,dy\,\, \times\int   t^4\,\,\vert K(t)\vert   \,\,dt \\
&\leqslant\frac{\beta_n^3n^{-7/2} k_n^4b_n^{-1}}{6}   \int_{\left\{f(y)\geqslant b_{n}+C\tau_{n}\right\}} \frac{ \left\vert R_{m,\ell}\left(y\right)R_{s,j}\left(y\right)\right\vert}{f^3(y)}\,\,dy\,\, \times\int   t^4\,\,\vert K(t)\vert   \,\,dt \\
&\leqslant\frac{\beta_n^3n^{-7/2} k_n^4b_n^{-1}}{6c_0^4}   \int  \left\vert R_{m,\ell}\left(y\right)R_{s,j}\left(y\right)\right\vert  f(y)\,\,dy\,\, \times\int   t^4\,\,\vert K(t)\vert   \,\,dt \\
&\leqslant \frac{n^{-7/2} k_n^4b_n^{-1}}{6c_0^4} \mathbb{E}\big[ \left\vert R_{m,\ell}\left(Y\right)R_{s,j}\left(Y\right)\right\vert\big]   \int   t^4\,\,\vert K(t)\vert   \,\,dt.
\end{align*}
Since $n^{-7/2} k_n^4b_n^{-1}\sim n^{-7/2+4c_{1}+c_{2}}=n^{-4(7/8-c_1-c_2/4)}\rightarrow 0$ as $n\rightarrow\infty$, it follows  from the preceding inequality that 
\begin{equation}\label{cv2l15}
\beta_n\sqrt{n}  \int\int_{\left\{f(y)\geqslant b_{n}+C\tau_{n}\right\}} R_{m,\ell}\left(y\right)R_{s,j}\left(y\right) \bigg(f\left(y+\frac{k_n\sqrt{\beta_n}}{nf(y)}t\right)-f(y)\bigg)\,K(t)\,\,dy\,dt=\mathrm{o}(1).
\end{equation}
Then, from \eqref{decl15}, \eqref{cv1l15} and \eqref{cv2l15}, we deduce that
\begin{equation*}
\sqrt{n}\mathbb{E}\left[\frac{R_{m,\ell}\left(Y\right)R_{s,j}\left(Y\right)\widehat{f}_{1,n}(Y)}{f(Y)}\mathbf{1}_{\left\{f(Y)\geqslant b_{n}+C\tau_{n}\right\}}\right]=\sqrt{n}\mathbb{E}\left[R_{m,\ell}\left(Y\right)R_{s,j}\left(Y\right)\right]+\mathrm{o}(1).
\end{equation*}
 Since $\widehat{f}_{2,n}$ is obtained from $\widehat{f}_{1,n}$ by replacing $\beta_n$ by $1/{\beta_n}$, we obtain the case of $r=2$  from a similar reasoning. It leads to 
\begin{align*}\label{dec2l15}
&\sqrt{n}\mathbb{E}\left[\frac{R_{m,\ell}\left(Y\right)R_{s,j}\left(Y\right)\widehat{f}_{2,n}(Y)}{f(Y)}\mathbf{1}_{\left\{f(Y)\geqslant b_{n}+C\tau_{n}\right\}}\right]\nonumber\\
&=  \frac{\sqrt{n}}{\beta_n}  \int\int_{\left\{f(y)\geqslant b_{n}+C\tau_{n}\right\}} R_{m,\ell}\left(y\right)R_{s,j}\left(y\right) \bigg(f\left(y+\frac{k_n}{n\sqrt{\beta_n}f(y)}t\right)-f(y)\bigg)\,K(t)\,\,dy\,dt\nonumber\\
&\hspace{0.5cm}+  \frac{\sqrt{n}}{\beta_n}  \int_{\left\{f(y)\geqslant b_{n}+C\tau_{n}\right\}} R_{m,\ell}\left(y\right)R_{s,j}\left(y\right)  f(y)\,\,dy,
\end{align*} 
with
\begin{align*}
&\bigg\vert\frac{\sqrt{n}}{\beta_n}  \int\int_{\left\{f(y)\geqslant b_{n}+C\tau_{n}\right\}} R_{m,\ell}\left(y\right)R_{s,j}\left(y\right) \bigg(f\left(y+\frac{k_n}{n\sqrt{\beta_n}f(y)}t\right)-f(y)\bigg)\,K(t)\,\,dy\,dt\bigg\vert\\
&\leqslant\frac{n^{-7/2} k_n^4b_n^{-1}}{6\beta_n^3c_0^4}   \int  \left\vert R_{m,\ell}\left(y\right)R_{s,j}\left(y\right)\right\vert  f(y)\,\,dy\,\, \times\int   t^4\,\,\vert K(t)\vert   \,\,dt \\
&\leqslant \frac{3n^{-7/2} k_n^4b_n^{-1}}{12c_0^4} \mathbb{E}\big[ \left\vert R_{m,\ell}\left(Y\right)R_{s,j}\left(Y\right)\right\vert\big]   \int   t^4\,\,\vert K(t)\vert   \,\,dt,
\end{align*}
for $n$ large enough, and
\begin{align*}
&\frac{\sqrt{n}}{\beta_n} \int_{\left\{f(y)\geqslant b_{n}+C\tau_{n}\right\}} R_{m,\ell}\left(y\right)R_{s,j}\left(y\right)  f(y)\,\,dy-\sqrt{n}\mathbb{E}\left[R_{m,\ell}\left(Y\right)R_{s,j}\left(Y\right)\right]\nonumber\\
&=\frac{(1-\beta_n)\sqrt{n}}{\beta_n}\mathbb{E}\left[R_{m,\ell}\left(Y\right)R_{s,j}\left(Y\right)\right]-\frac{1}{\beta_n}\sqrt{n}\mathbb{E}\left[R_{m,\ell}\left(Y\right)R_{s,j}\left(Y\right)\mathbf{1}_{\left\{f(Y)< b_{n}+C\tau_{n}\right\}}\right]\\
&=\mathrm{o}(1).
\end{align*}
Hence
\begin{equation*}
\sqrt{n}\mathbb{E}\left[\frac{R_{m,\ell}\left(Y\right)R_{s,j}\left(Y\right)\widehat{f}_{2,n}(Y)}{f(Y)}\mathbf{1}_{\left\{f(Y)\geqslant b_{n}+C\tau_{n}\right\}}\right]=\sqrt{n}\mathbb{E}\left[R_{m,\ell}\left(Y\right)R_{s,j}\left(Y\right)\right]+\mathrm{o}(1).
\end{equation*}
\hfill $\Box$

\bigskip

\begin{Lem}\label{l16}
Under the assumptions \ref{as2}, \ref{as9}, \ref{as5} and \ref{as6},   we have for  $(\ell,j)\in\{1,\ldots,d\}^2$:
\begin{equation*}
\sqrt{n}\mathbb{E}\left[\frac{R_{b_n,\ell}\left(Y\right)R_{b_n,j}\left(Y\right)\widehat{f}_{b_{n}}(Y)}{f_{b_{n}}(Y)}\right]=\sqrt{n}\mathbb{E}\left[R_{\ell}\left(Y\right)R_{j}\left(Y\right)\right]+\mathrm{o}(1).
\end{equation*}
\end{Lem}
\noindent\textit{Proof.} Arguing as in \cite{zhu96} (see Eq. (4.21), p. 1066), we have 
\[
\widehat{f}_{b_{n}}(Y)=J_{1,n}(Y)+J_{2,n}(Y)+J_{3,n}(Y)+J_{4,n}(Y)+J_{5,n}(Y),
\]
where
\begin{align*}
J_{1,n}(Y)&=\widehat{f}_{n}(Y)\mathbf{1}_{\left\{f(Y)\geqslant b_{n}+C\tau_{n}\right\}},\notag\\
J_{2,n}(Y)&=f(Y)\left\{\mathbf{1}_{\left\{\widehat{f}_{n}(Y)\geqslant b_{n}\right\}}-\mathbf{1}_{\left\{f(Y)\geqslant b_{n}+C\tau_{n}\right\}}\right\},\nonumber\\
J_{3,n}(Y)&=\left(\widehat{f}_{n}(Y)-f(Y)\right)\left\{\mathbf{1}_{\left\{\widehat{f}_{n}(Y)\geqslant b_{n}\right\}}-\mathbf{1}_{\left\{f(Y)\geqslant b_{n}+C\tau_{n}\right\}}\right\},\nonumber\\
J_{4,n}(Y)&=b_{n}\mathbf{1}_{\left\{f(Y)< b_{n}-C\tau_{n}\right\}},\notag\\
J_{5,n}(Y)&=b_{n}\left\{\mathbf{1}_{\left\{\widehat{f}_{n}(Y)< b_{n}\right\}}-\mathbf{1}_{\left\{f(Y)< b_{n}+C\tau_{n}\right\}}\right\},\notag
\end{align*}
$C$ is a given positive constant and $\tau_n= \frac{k_{n}^4}{n^4}+\frac{\sqrt{n\log(n)}}{k_{n}}$. Then, it is enough to show that 
\begin{equation}\label{j1n}
\sqrt{n}\mathbb{E}\left[\frac{R_{b_n,\ell}\left(Y\right)R_{b_n,j}\left(Y\right)J_{1,n}(Y)}{f_{b_{n}}(Y)}\right]=\sqrt{n}\mathbb{E}\left[R_{\ell}\left(Y\right)R_{j}\left(Y\right)\right]+\mathrm{o}(1)
\end{equation}
and
\begin{equation}\label{jmn}
\sqrt{n}\mathbb{E}\left[\frac{R_{b_n,\ell}\left(Y\right)R_{b_n,j}\left(Y\right)J_{m,n}(Y)}{f_{b_{n}}(Y)}\right]=\mathrm{o}(1), \,\,m=2,\ldots,5.
\end{equation}
Proof of \eqref{j1n}: Clearly,
\begin{align}\label{deczet}
&\sqrt{n}\mathbb{E}\left[\frac{R_{b_n,\ell}\left(Y\right)R_{b_n,j}\left(Y\right)J_{1,n}(Y)}{f_{b_{n}}(Y)}\right]\nonumber\\
&=\sqrt{n}\mathbb{E}\left[\frac{R_{\ell}\left(Y\right)R_{j}\left(Y\right)\widehat{f}_{n}(Y)}{f(Y)}\mathbf{1}_{\left\{f(Y)\geqslant b_{n}+C\tau_n \right\}}\right]\nonumber\\
&=\sqrt{n}\mathbb{E}\left[\frac{(R_{1,\ell}(Y)-R_{2,\ell}(Y))(R_{1,j}(Y)-R_{2,j}(Y))\widehat{f}_{n}(Y)}{f(Y)}\mathbf{1}_{\left\{f(Y)\geqslant b_{n}+C\tau_n \right\}}\right]\nonumber\\
&=\zeta_{1,1,n}-\zeta_{1,2,n}-\zeta_{2,1,n}+\zeta_{2,2,n},
\end{align}
where 
\[
\zeta_{m,s,n}=\sqrt{n}\mathbb{E}\left[\frac{R_{m,\ell}(Y) R_{s,j}(Y) \widehat{f}_{n}(Y)}{f(Y)}\mathbf{1}_{\left\{f(Y)\geqslant b_{n}+C\tau_n \right\}}\right].
\]
Since, from \cite{bengono24} (see p. 12),  $\widehat{f}_{1,n}\leqslant \widehat{f}_{n}\leqslant\widehat{f}_{2,n}$,  where $\widehat{f}_{1,n}$ and $\widehat{f}_{2,n}$ are defined in \eqref{f1n} and \eqref{f2n}, we deduce that
\begin{align*}
\sqrt{n}\mathbb{E}\left[\frac{R_{m,\ell}(Y) R_{s,j}(Y) \widehat{f}_{1,n}(Y)}{f(Y)}\mathbf{1}_{\left\{f(Y)\geqslant b_{n}+C\tau_n \right\}}\right]
\leqslant \zeta_{m,s,n}  \leqslant \sqrt{n}\mathbb{E}\left[\frac{R_{m,\ell}(Y) R_{s,j}(Y) \widehat{f}_{2,n}(Y)}{f(Y)}\mathbf{1}_{\left\{f(Y)\geqslant b_{n}+C\tau_n \right\}}\right],
\end{align*}
and  from Lemma \ref{l15} that
\begin{align*}
\mathrm{o}(1)
\leqslant \zeta_{m,s,n}-\sqrt{n}\mathbb{E}\left[R_{m,\ell}\left(Y\right)R_{s,j}\left(Y\right)\right]  \leqslant \mathrm{o}(1),
\end{align*}
what implies that
\begin{equation}\label{zetms}
\zeta_{m,s,n}=\sqrt{n}\mathbb{E}\left[R_{m,\ell}\left(Y\right)R_{s,j}\left(Y\right)\right]+\mathrm{o}(1).
\end{equation}
Then, from \eqref{deczet}, \eqref{zetms} and \eqref{decR}, we obtain \eqref{j1n}.

\bigskip

\noindent Proof of \eqref{jmn}:  Arguing as in \cite{zhu96} (see pp. 1066--1067), and using the fact that $b_{n}+C\tau_n \sim b_n$ since $b_n^{-1}\tau_n\rightarrow 0$ as $n\rightarrow\infty$, we get the inequalities
\[
\bigg\vert\sqrt{n}\mathbb{E}\left[\frac{R_{b_n,\ell}\left(Y\right)R_{b_n,j}\left(Y\right)J_{2,n}(Y)}{f_{b_{n}}(Y)}\right]\bigg\vert\leqslant \sqrt{n}\mathbb{E}\left[\left|R_{\ell}(Y)R_{j}(Y)\right|\mathbf{1}_{\left\{f(Y)< b_{n}+C\tau_n\right\}}\right]=\mathrm{o}(1),
\]
\[
\bigg\vert\sqrt{n}\mathbb{E}\left[\frac{R_{b_n,\ell}\left(Y\right)R_{b_n,j}\left(Y\right)J_{3,n}(Y)}{f_{b_{n}}(Y)}\right]\bigg\vert\leqslant C\tau_n b_{n}^{-1}\sqrt{n}\mathbb{E}\left[\left|R_{\ell}(Y)R_{j}(Y)\right|\mathbf{1}_{\left\{f(Y)< b_{n}+C\tau_n\right\}}\right] =\mathrm{o}(1),
\]
\[
\bigg\vert\sqrt{n}\mathbb{E}\left[\frac{R_{b_n,\ell}\left(Y\right)R_{b_n,j}\left(Y\right)J_{4,n}(Y)}{f_{b_{n}}(Y)}\right]\bigg\vert\leqslant \bigg\vert\sqrt{n}\mathbb{E}\left[R_{\ell}(Y)R_{j}(Y)\mathbf{1}_{\left\{f(Y)< b_{n}-C\tau_n\right\}}\right]\bigg\vert=\mathrm{o}(1),
\]
\[
\bigg\vert\sqrt{n}\mathbb{E}\left[\frac{R_{b_n,\ell}\left(Y\right)R_{b_n,j}\left(Y\right)J_{5,n}(Y)}{f_{b_{n}}(Y)}\right]\bigg\vert\leqslant  \bigg\vert\sqrt{n}\mathbb{E}\left[R_{\ell}(Y)R_{j}(Y)\mathbf{1}_{\left\{f(Y)< b_{n}+C\tau_n\right\}}\right] \bigg\vert=\mathrm{o}(1), 
\]
which complete the proof of  \eqref{jmn}.
\hfill $\Box$

\subsection{Proof of Theorem \ref{theo1}}
Let us denote by  $\widehat{\lambda}_{\ell j}^{(n)}$ the $(\ell,j)$-th entry  of the  $d\times d$  matrix $\widehat{\Lambda}_n$.  It is easily seen that
\begin{eqnarray*}
\sqrt{n}\,\widehat{\lambda}_{\ell j}^{(n)}=\frac{1}{\sqrt{n}}\sum_{i=1}^n\frac{\widehat{g}_{\ell , n}(Y_i)\widehat{g}_{j , n}(Y_i)}{\widehat{f}_{b_n}^2(Y_i)}&=&\frac{1}{\sqrt{n}}\sum_{i=1}^n\Big\{I_{\ell j}^{(1)}(Y_i)+I_{\ell j}^{(2)}(Y_i)-I_{\ell j}^{(3)}(Y_i)\Big\} \\
&&- U_{n,\ell,j}^{(1)}+U_{n,\ell,j}^{(2)} + U_{n,\ell,j}^{(3)} - U_{n,\ell,j}^{(4)},
\end{eqnarray*}
where
\begin{equation*} 
I_{\ell j}^{(1)}(y) =\frac{g_\ell(y)g_j(y)}{f_{b_n}^2(y)}= R_{b_n,\ell}(y)R_{b_n,j}(y),
\end{equation*}
$I_{\ell j}^{(2)}$ and $I_{\ell j}^{(3)}$ are defined in \eqref{Ilj2} and \eqref{Ilj3}, and the $U_{n,\ell,j}^{(m)}$'s  are given in \eqref{u1}, \eqref{u2}, \eqref{u3} and \eqref{u4}.
Then,  from Lemma  \ref{cl14}  we get
\begin{equation}\label{b22} 
  \sqrt{n}\,\widehat{\lambda}_{\ell j}^{(n)}=\frac{1}{\sqrt{n}}\sum_{i=1}^n\Big\{I_{\ell j}^{(1)}(Y_i)+I_{\ell j}^{(2)}(Y_i)-I_{\ell j}^{(3)}(Y_i)\Big\} +\mathrm{o}_p(1),
\end{equation}
and, putting $\nu_{\ell j}=\mathbb{E}\left(I_{\ell j}^{(1)}(Y)+I_{\ell j}^{(2)}(Y)-I_{\ell j}^{(3)}(Y)\right)$, we have 
\begin{align*}
&\sqrt{n}\left(\widehat{\lambda}_{\ell j}^{(n)}-\nu_{\ell j}\right)\\
&=\frac{1}{\sqrt{n}}\sum_{i=1}^n\Big\{I_{\ell j}^{(1)}(Y_i)-\mathbb{E}\left(I_{\ell j}^{(1)}(Y)\right)\Big\}+\frac{1}{\sqrt{n}}\sum_{i=1}^n\Big\{I_{\ell j}^{(2)}(Y_i)-\mathbb{E}\left(I_{\ell j}^{(2)}(Y)\right)\Big\}\\
&-\frac{1}{\sqrt{n}}\sum_{i=1}^n\Big\{I_{\ell j}^{(3)}(Y_i)-\mathbb{E}\left(I_{\ell j}^{(3)}(Y)\right)\Big\}+\mathrm{o}_p(1). 
\end{align*}
Then, from   Lemma \ref{l12to13} and the three following  equality obtained in Eq. (4.6), Eq. (4.14) and Eq. (4.15) of  \cite{zhu96} :
\[
\frac{1}{\sqrt{n}}\sum_{i=1}^n\Big\{I_{\ell j}^{(1)}(Y_i)-\mathbb{E}\left(I_{\ell j}^{(1)}(Y)\right)\Big\}=\frac{1}{\sqrt{n}}\sum_{i=1}^n\Big\{R_\ell(Y_i)R_j(Y_i)-\mathbb{E}\left[ R(Y_\ell)R_j(Y)\right]\Big\}+\mathrm{o}_p(1),
\]
\begin{align*}
&\frac{1}{\sqrt{n}}\sum_{i=1}^n\bigg\{ R_{b_n,\ell}(Y_i)R_{b_n,j}(Y_i)+\frac{1}{2}X_{i\ell}R_{b_n,j}(Y_i)\frac{f(Y_i)}{f_{b_n}(Y_i)}+\frac{1}{2}X_{ij}R_{b_n,\ell}(Y_i)\frac{f(Y_i)}{f_{b_n}(Y_i)}\\
&-\mathbb{E}\left( R_{b_n,\ell}(Y)R_{b_n,j}(Y)+\frac{1}{2}X_{\ell}R_{b_n,j}(Y)\frac{f(Y)}{f_{b_n}(Y)}+\frac{1}{2}X_{j}R_{b_n,\ell}(Y)\frac{f(Y)}{f_{b_n}(Y)}\right)\bigg\}\\
&=\frac{1}{\sqrt{n}}\sum_{i=1}^n\Big\{R_\ell(Y_i)R_j(Y_i)+\frac{1}{2}X_{i\ell}R_{j}(Y_i)+\frac{1}{2}X_{ij}R_{\ell}(Y_i)-2\mathbb{E}\left[ R(Y_\ell)R_j(Y)\right]\Big\}+\mathrm{o}_p(1),
\end{align*}
and
\begin{align*}
&\frac{1}{\sqrt{n}}\sum_{i=1}^n\bigg\{ R_{b_n,\ell}(Y_i)R_{b_n,j}(Y_i)\frac{f(Y_i)}{f_{b_n}(Y_i)}-\mathbb{E}\left( R_{b_n,\ell}(Y)R_{b_n,j}(Y)\frac{f(Y)}{f_{b_n}(Y)}\right)\bigg\}\\
&=\frac{1}{\sqrt{n}}\sum_{i=1}^n\Big\{R_\ell(Y_i)R_j(Y_i)-\mathbb{E}\left[ R_\ell(Y)R_j(Y)\right]\Big\}+\mathrm{o}_p(1),
\end{align*}
it follows
\[
\sqrt{n}\left(\widehat{\lambda}_{\ell j}^{(n)}-\nu_{\ell j}\right)=\frac{1}{\sqrt{n}}\sum_{i=1}^n\bigg\{\frac{1}{2}\left(X_{i\ell}R_j(Y_i)+X_{ij}R_{\ell}(Y_i)\right)-\mathbb{E} \left(R_\ell(Y)R_{j}(Y)\right)\bigg\}+\mathrm{o}_p(1).
\]
Moreover,
\[
\sqrt{n}\,\nu_{k,\ell}=\sqrt{n}\mathbb{E}\bigg[  R_{b_n,\ell}(Y)R_{b_n,j}(Y)+ R_{b_n,\ell}(Y)\frac{\widehat{g}_{j,n}(Y)}{f_{b_n}(Y)}+ R_{b_n,j}(Y)\frac{\widehat{g}_{\ell,n}(Y)}{f_{b_n}(Y)}-2R_{b_n,\ell}(Y)R_{b_n,j}(Y)\frac{\widehat{f}_{b_n}(Y)}{f_{b_n}(Y)}\bigg];
\]
then, using  Lemma \ref{l15}, Lemma \ref{l16} together with the equality
\[
\sqrt{n}\mathbb{E}\bigg[  R_{b_n,\ell}(Y)R_{b_n,j}(Y)\bigg]=\sqrt{n}\mathbb{E}\left[ R(Y_\ell)R_j(Y)\right]+\mathrm{o}(1)
\]
given in Eq. (4.17) of \cite{zhu96}, we obtain $\sqrt{n}\,\nu_{\ell j}=\sqrt{n}\,\lambda_{\ell j}+\mathrm{o}(1)$, where $\lambda_{\ell j}=\mathbb{E} \left(R_\ell(Y)R_{j}(Y)\right)$. Hence
\[
\sqrt{n}\left(\widehat{\lambda}_{\ell j}^{(n)}-\lambda_{\ell j}\right)=\frac{1}{\sqrt{n}}\sum_{i=1}^n\bigg\{\frac{1}{2}\left(X_{i\ell }R_j(Y_i)+X_{ij}R_{\ell}(Y_i)\right)-\mathbb{E} \left(R_\ell(Y)R_{j}(Y)\right)\bigg\}+\mathrm{o}_p(1).
\]
Clearly,
\[
\mathbb{E} \left(\frac{1}{2}\left(X_{\ell }R_j(Y)+X_{j}R_{\ell}(Y)\right)\right)=\mathbb{E} \Big(X_{\ell}R_j(Y)\Big) =\mathbb{E} \Big(R_{\ell}(Y)R_j(Y)\Big),
\]
and, putting  $\mathcal{H}_{n}=\sqrt{n}\,\left(\widehat{\Lambda}_n-\Lambda\right)$ and $\mathcal{H}^{(n)}_{\ell j}=\sqrt{n}\left(\widehat{\lambda}_{\ell j}^{(n)}-\lambda_{\ell j}\right)$, we have
\[
\textrm{Tr}\left(A^\top\mathcal{H}_{n}\right)=\sum_{\ell=1}^d\sum_{j=1}^da_{\ell j}\,\mathcal{H}^{(n)}_{\ell j}
=\frac{1}{\sqrt{n}}\sum_{i=1}^n\left(\mathcal{U}_i-\mathbb{E} \left(\mathcal{U}_i\right)\right)+o_p(1),
\]
where 
\[
\mathcal{U}_i=\sum_{\ell=1}^d\sum_{j=1}^d\frac{a_{\ell j}}{2}\left(X_{i\ell }R_j(Y_i)+X_{ij}R_{\ell}(Y_i)\right).
\]
From the central limit theorem and  Slutsky's theorem we deduce  that $\textrm{Tr}\left(A^\top\mathcal{H}_{n}\right)\stackrel{\mathscr{D}}{\rightarrow}\mathcal{N}\left(0,\sigma_{A}^{2}\right)$,   as $n\rightarrow \infty$, where 
$
\sigma_A^2  
$ 
is given in \eqref{siga}. Then, using Levy's theorem, we conclude that  $ \mathcal{H}_{n}\stackrel{\mathscr{D}}{\rightarrow}\mathcal{H}$, as $n\rightarrow \infty$, where $\mathcal{H}$ has a normal distribution in $\mathscr{M}_d(\mathbb{R})$ with   $\textrm{Tr}\left(A^\top\mathcal{H}_{n}\right)\leadsto  \mathcal{N}\left(0,\sigma_{A}^{2}\right)$.


\begin{thebibliography}{9}



\bibitem{bengono24}
L. Bengono Mintogo, E.D.D. Nkou, and G.M. Nkiet,
Rates of strong uniform consistency for the $k$-nearest neighbors kernel estimators of density and regression function, ArXiv:2408.12741 (2024).

\bibitem{bura01}
E. Bura and    R.D.  Cook,  Estimating the structural dimension of regressions via parametric inverse regression, J. R. Stat. Soc. Ser. B. Stat. Methodol. \textbf{63}(2001), 393--410.


\bibitem{collomb79}
G. Collomb, Estimation de la r\'egression par la m\'ethode des $ k$  points les plus proches avec noyau: Quelques
propriétés de convergence ponctuelle,  C. R. Math. Acad. Sci. Paris  \textbf{289}(1979), 245--247.



\bibitem{cook02}
R.D. Cook and B. Li,  
Dimension reduction for conditional mean in regression,  Ann.  Statist. \textbf{30}(2002),  455--474.

\bibitem{cook91}
R.D. Cook and S. Weisberg,  
Comment on ‘Sliced inverse regression for dimension reduction’, J. Amer. Statist. Assoc.  \textbf{86}(1991),  328--332.

\bibitem{coudret14}
R. Coudret, S. Girard, and J. Sarraco,
A new sliced inverse regression method for multivariate response,  Comput. Statist. Data Anal.  \textbf{77}(2014), 285--299.



\bibitem{ebende24}
J.R. Ebende Penda, E.D.D. Nkou, S. Bouka, and G.M. Nkiet, 
On variable selection in partially linear regression model,  C. R. Math. Acad. Sci. Soc. R. Can. \textbf{46} (2024), 119--144.

\bibitem{ferre98}
L. Ferr\'e, 
Determining the dimension in sliced inverse regression and related methods,  J. Amer. Statist. Assoc. \textbf{93}(1998),  132--140.



\bibitem{hsing92}
T. Hsing and R.J. Caroll,
An asymptotic theory for sliced inverse regression,   Ann.  Statist. \textbf{20}(1992),  1040--1061.

\bibitem{li91}
K.C. Li,  
Sliced inverse regression for dimension reduction,  J. Amer. Statist. Assoc.  \textbf{86}(1991),  316--342.

\bibitem{li03}
K.C. Li,  Y. Aragon, K. Shedden, and C. Thomas-Agnan,
Dimension reduction for multivariate response data,  J. Amer. Statist. Assoc.  \textbf{98}(2003),  99--109.

\bibitem{moore77}
D. S. Moore and J. W. Yackel, 
Consistency properties of nearest neighbor density function estimators, Ann. Statist. \textbf{5}(1977),  143--154.

\bibitem{nkiet08}
G.M. Nkiet,  
Consistent estimation of the dimensionality in sliced inverse regression,  Ann. Inst. Statist. Math.   \textbf{60}(2008), 257--271.

\bibitem{nkou23}
E.D.D. Nkou, 
Recursive kernel estimator in semiparametric regression model,  J. Nonparametr. Stat.  \textbf{35}(2023),   145--171.

\bibitem{nkou19}
E.D.D. Nkou and G.M. Nkiet,  
Strong consistency of kernel estimator in a semiparametric regression model,  Statistics \textbf{53}(2019), 1289--1305.

\bibitem{nkou22}
E.D.D. Nkou and G.M. Nkiet,  
Wavelet-based estimation   in a semiparametric regression model,  Int. J. Wavelets Multiresolut. Inf. Process. \textbf{20}(2022), 2150056.


\bibitem{powell89}
J.L. Powell, J.H.  Stock,  and T.M. Stoker,
Semiparametric estimation of index coefficients,   Econometrica  \textbf{57}(1989), 1403--1430.

\bibitem{schott94}
J.R. Schott,   
Determining the dimensionality in sliced inverse regression, J. Amer. Statist. Assoc.  \textbf{89}(1994), 141--148.

\bibitem{velilla98}
S. Velilla, 
Assessing the number of linear components in a general regression problem.
J. Amer. Statist. Assoc.  \textbf{93}(1998),   1088--1098.

\bibitem{zhu96}
 L.-X.  Zhu and   K.-T.  Fang, 
Asymptotics for kernel estimate of sliced inverse regression,  Ann. Statist. \textbf{24}(1996), 1053--1068.


\bibitem{zhu06}
L. Zhu,  B. Miao,   and H. Peng,
On sliced inverse regression with high dimensional covariates, J. Amer. Statist. Assoc.  \textbf{101}(2006),  630--643.

\bibitem{zhu95}
 L.-X.  Zhu and   K.W. Ng, 
Asymptotics  of sliced inverse regression, Statist.  Sinica \textbf{5}(1995), 727--736.

\bibitem{zhu07}
L.-P. Zhu and L.-X. Zhu,  On kernel method for sliced average variance
estimation,  J. Multivariate Anal. \textbf{98}(2007), 970–991.


\end{thebibliography}
\end{document}